\documentclass[reqno,psamsfonts, 10pt]{amsart}
\usepackage{graphicx, amssymb, amsmath, amsfonts, hyperref, amsthm, amsbsy}
\usepackage{mathrsfs}

\newtheorem{thm}{Theorem}[section]
\newtheorem{cor}[thm]{Corollary}
\newtheorem{lem}[thm]{Lemma}

\theoremstyle{remark}
\newtheorem{rem}[thm]{Remark}

\newcommand{\cov}{\mbox{$\mathrm{cov}$}}

\newcommand{\sinc}{\mbox{$\mathrm{sinc}$}}

\newcommand{\comment}[1]{}

\numberwithin{equation}{section}

\begin{document}

\date{\today}

\title[A new REM conjecture]{A new REM conjecture }
\author[G. Ben Arous]{G. Ben Arous}
\address{G. Ben Arous\\Swiss Federal Institute of Technology
(EPFL), CH-1015 Lausanne, Switzerland and Courant Institute of
Mathematical Sciences, New York University, 251 Mercer Street, New
York NY 10012} \email{benarous@cims.nyu.edu}

\author[V. Gayrard]{V\'eronique Gayrard}
\address{V\'eronique Gayrard\\Laboratoire d'Analyse, Topologie, Probabilit\'es\\
CMI, 39 rue Joliot-Curie\\
13453 Marseille Cedex} \email{gayrard@latp.univ-mrs.fr,
veronique@gayrard.net}

\author[A. Kuptsov]{A. Kuptsov}
\address{A. Kuptsov\\Swiss Federal Institute of Technology
(EPFL), CH-1015 Lausanne, Switzerland and Courant Institute of
Mathematical Sciences, New York University, 251 Mercer Street, New
York NY 10012} \email{kuptsov@cims.nyu.edu}


\begin{abstract}
We introduce here a new universality conjecture for levels of random
Hamiltonians, in the same spirit as the local REM conjecture made by
S. Mertens and H. Bauke. We establish our conjecture for a wide
class of Gaussian and non-Gaussian Hamiltonians, which include the
$p$-spin models, the Sherrington-Kirkpatrick model and the number
partitioning problem. We prove that our universality result is
optimal for the last two models by showing when this universality
breaks down.
\end{abstract}

\maketitle

\section{Introduction}

S. Mertens and H. Bauke recently observed (\cite{Mer00},
\cite{BM04}, see also \cite{BFM04}) that the statistics of energy
levels for very general random Hamiltonians are Poissonian, when
observed micro-canonically, i.e. in a small window in the bulk. They
are universal and identical to those of the simplest spin-glass
model, the Random Energy Model or REM, hence the name of this
(numerical) observation: the REM conjecture or more precisely the
\emph{local REM conjecture}.

This local REM conjecture was made for a wide class of random
Hamiltonians of statistical mechanics of disordered systems, mainly
spin-glasses (mean field or not), and for various combinatorial
optimization problems (like number partitioning). Recently, two
groups of mathematicians have established this conjecture in
different contexts : C. Borgs, J. Chayes, C.Nair, S.Mertens,
B.Pittel for the number partitioning question (see \cite{BCP},
\cite{BCMN1},\cite{BCMN2}), and A. Bovier, I. Kurkova for general
spin-glass Hamiltonians (see \cite{BK}, \cite{BK2}).

We introduce here a new kind of universality for the energy levels
of disordered systems. We believe that one should find universal
statistics for the energy levels of a wide class of random
Hamiltonians if one re-samples the energy levels, i.e. draws a
random subset of these energies. Put otherwise, our conjecture is
thus that the level statistics should also be universal, i.e
Poissonian, when observed on a large random subset of the
configuration space rather than in a micro-canonical window.  We
establish this new universality result (which could be called
\emph{the re-sampling REM universality} or \emph{the REM
universality by dilution}) for general (mean-field) spin-glass
models, including the case of number partitioning and for large but
sparse enough subsets. This approach has the following interesting
property: the range of energies involved is not reduced to a small
window as in the local REM conjecture. Thus we can study the extreme
value distribution on the random subset, by normalizing the energies
properly. Doing so we establish that the Gibbs measure restricted to
a sparse enough random subset of configuration space has a universal
distribution which is thus the same as for the REM, i.e. a
Poisson-Dirichlet measure.

To be more specific, we specialize our setting to the case of random
Hamiltonians $(H_N(\sigma))_{\sigma \in S_N}$ defined on the
hypercube $S_N=\{-1,1\}^N$.  We want to consider a sparse random
subset of the hypercube, say $X$, and the restriction of the
function $H_N$ to $X$. We introduce the random point process
\begin{equation}
\mathcal{P}_N = \sum \limits_{\sigma \in X} \delta_{H'_N(\sigma)},
\end{equation}
with a normalization:
\begin{equation}\label{(intro):Hamiltonian:normalized}
H'_N(\sigma)=\frac{H_N(\sigma)-a_N}{b_N}
\end{equation}
to be chosen. Our conjecture specializes to the following: the
asymptotic behavior of the random point process $\mathcal{P}_N$ is
universal, for a large class of Hamiltonians $H_N$, for appropriate
sparse random subsets $X$, and appropriate normalization.

We will only study here the simplest possible random subset, i.e.
a site percolation cluster $X=\{\sigma \in S_N: X_{\sigma}=1\}$,
where the random variables $(X_{\sigma})_{\sigma \in S_N}$ are
i.i.d. and Bernoulli:
\begin{equation}
\mathrm{P}(X_{\sigma}=1)=1-\mathrm{P}(X_{\sigma}=0)=: p_{N}
=\frac{2^M}{2^N}.
\end{equation}
Thus the mean size of $X$ is $2^M$ and we will always assume that
$X$ is not too small, e.g. that $\log N = o(2^M).$ We will sometimes
call $X$ a \textit{random cloud}.

In order to understand what the universal behavior should be, let us
examine the trivial case where the $H_N(\sigma)$ are i.i.d centered
standard Gaussian random variables, i.e. the case of the Random
Energy Model. Then, standard extreme value theory proves that if
\begin{align}\label{normalisation}
a_N = \sqrt{2M \log 2 + 2 \log b_N - \log 2 } \mbox{ and } b_N =
\sqrt{\frac{1}{M}},
\end{align}
then $\mathcal{P}_N$ converges to a Poisson point process with
intensity measure
\begin{equation}
\mu(dt) = \frac{1}{\sqrt{\pi}} e^{-t \sqrt{2 \log 2}} dt.
\end{equation}
We will now fix the normalization needed in
\eqref{(intro):Hamiltonian:normalized} by choosing $a_N$ and $b_N$ as
in \eqref{normalisation}. The basic mechanism of the REM
universality we propose is that the influence of correlations
between the random variables $H'_N(\sigma)$ should be negligible
when the two-point correlation (the covariance) is a decreasing
function of the Hamming distance
\begin{equation}
d_H(\sigma,\sigma') = \# \{ i \le N: \sigma_i \neq \sigma'_i \}
\end{equation}
and when the random cloud is sparse enough. We first establish this
universality conjecture  for a large class of Gaussian Hamiltonians.
This class contains the Sherrington-Kirkpatrick (SK) model as well
as the more general $p$-spin models. It also contains the Gaussian
version of the number partitioning problem.

\par Consider a Gaussian Hamiltonian
$H_N$ on the hypercube $S_N=\{-1,1\}^N$ such that the random
variables $(H_N(\sigma))_{\sigma \in S_N}$ are centered and whose
covariance is a smooth decreasing function  of the Hamming distance
or, equivalently a smooth increasing function $\nu$ of the overlap
$R(\sigma,\sigma')= \frac{1}{N} \sum \limits_{i=1}^{N} \sigma_i
\sigma'_i $ :
\begin{equation}\label{(intro):cov:str}
\cov(H_N(\sigma),H_N(\sigma')) = \nu(R(\sigma,\sigma'))= \nu
\Bigl(1 - \frac{2 d_H(\sigma,\sigma')}{N} \Bigr).
\end{equation}
We will always assume that $\nu(0)=0$ and that $\nu(1)=1$. The first
assumption is crucial, since it means that the correlation of the
Hamiltonian vanishes for pairs of points on the hypercube which are
at a typical distance. The second assumption simply normalizes the
variance of $H_N(\sigma)$ to 1.

This type of covariance structure can easily be realized for $\nu$
real analytic, of the form:
\begin{equation}\label{(intro):cov:str'}
\nu(r)= \sum \limits_{p \geq 1} a_p^2 r^p.
\end{equation}

Indeed such a covariance structure can be realized by taking
mixtures of $p$-spin models. Let $H_{N,p}$ be the Hamiltonian of the
$p$-spin model given by
\begin{equation}\label{(intro):p-spin_SK}
H_{N,p}(\sigma) = \frac{1}{N^{{(p-1)}/{2}}} \sum \limits_{1 \leq
i_1, i_2, \dots, i_p \leq N} g_{i_1,i_2,\dots,i_p}\sigma_{i_1}
\sigma_{i_2}\dots\sigma_{i_p},
\end{equation}
where random variables $(g_{i_1,i_2,\dots,i_p})_{1 \leq
i_1,i_2,\dots,i_p \leq N}$ are independent standard Gaussians
defined on a common probability space $(\Omega^{g}, \mathcal{F}^{g},
\mathbb{P}).$ Then
\begin{equation}\label{(intro):Hamiltonian}
H_{N}(\sigma) = \frac{1}{\sqrt{N}} \sum \limits_{p \geq 1} a_p
H_{N,p}(\sigma)
\end{equation}
has the covariance structure given in
\eqref{(intro):cov:str}-\eqref{(intro):cov:str'} Let
us recall that the case where $\nu(r)=r$  is the Gaussian version of
the number partitioning problem (\cite{BCP}), the case where
$\nu(r)=r^2$  is the SK model, and more generally, $\nu(r)=r^p$
defines the pure $p$-spin model.

\noindent Let us normalize $H_N$ as above (see
\eqref{normalisation}):

\begin{equation}
H'_N(\sigma)=\frac{H_N(\sigma)-a_N}{b_N},
\end{equation}
and consider the sequence of point processes
\begin{equation}\label{1}
\mathcal{P}_N = \sum \limits_{\sigma \in X} \delta_{H'_N(\sigma)}.
\end{equation}

\begin{thm}[Universality in the Gaussian case]\label{main_thm}
Assume that $M = o(\sqrt{N})$ if $\nu'(0) \neq 0 $, and that $M =
o(N)$ if $\nu'(0)=0$. Then, $\mathrm{P}$-almost surely, the
distribution of the point process $\mathcal{P}_N$ converges weakly
to the distribution of a Poisson point process $\mathcal{P}$ on
$\mathbb{R}$ with intensity given by
\begin{equation}\label{(intro):intensity}
\mu(dt) = \frac{1}{\sqrt{\pi}} e^{-t \sqrt{2 \log 2}} dt.
\end{equation}
\end{thm}

\begin{rem}
The condition $\log N = o(2^M)$ is needed in order to get
$\mathrm{P}$-almost sure results.
\end{rem}

We extend this result, in Section \ref{NG}, to a wide class of
non-Gaussian Hamiltonians (introduced in \cite{BCMN1} for the case
of number partitioning).

\noindent The theorem has the following immediate corollary. Let us
fix the realization of the random cloud $X.$ For configurations
$\sigma$ belonging to the cloud we consider the Gibbs' weights
$G_{N,\beta}(\sigma)$ of the re-scaled Hamiltonian $H'_N(\sigma)$
\begin{equation}
G_{N,\beta}(\sigma) = \frac{e^{-\beta H'_N(\sigma)}}{\sum
\limits_{\varrho \in X} e^{-\beta H'_N(\varrho)}} = \frac{e^{-\beta
\sqrt{M} H_N(\sigma)}}{\sum \limits_{\varrho \in X} e^{-\beta
\sqrt{M} H_N(\varrho)}}.
\end{equation}
Reordering the Gibbs' weights $G_{N,\beta}(\sigma)$ of the
configurations $\sigma \in X$ as a non-increasing sequence
$(w_\alpha)_{\alpha \leq |X|}$ and defining $w_\alpha=0$ for $\alpha
> |X|$ we get a random element $w$ of the space $\mathcal{S}$
of non increasing sequences of non negative real numbers with sum
less than one.

\begin{cor}[Convergence to Poisson-Dirichlet]\label{(a):cor_Gibbs}
If $\beta > \sqrt{2 \log 2}$ then $\mathrm{P}$-almost surely under
the assumptions of Theorem \ref{main_thm} the law of the sequence
$w=(w_\alpha)_{\alpha \geq 1}$ converges to the Poisson-Dirichlet
distribution with parameter $m=\frac{\sqrt{\smash[b]2 \log
2}}{\beta}$ on $\mathcal{S}.$
\end{cor}
The fact that Theorem \ref{main_thm} implies Corollary
\ref{(a):cor_Gibbs} is well-known, see for instance \cite{Tal03}
(pp.13-19) for a good exposition.

\comment{ \par To prove the convergence in Theorem \ref{main_thm} we
use the classical approach known as the moment method. This method
was used in works of C. Borgs, J. Chayes, S. Mertens and C. Nair
\cite{BCMN1},~\cite{BCMN2} and A. Bovier, I. Kurkova \cite{BK},
\cite{BK2}. It is essentially based on the following fact. Let
$\mathscr{B}$ consist of all bounded Borel sets on the real line.
Suppose that for every set $A \in \mathscr{B}$ and for every $\ell
\in \mathbb{N}$ the $\ell^{\textrm{th}}$ factorial moment of random
variable \begin{equation} \mathcal{P}_N(A) = |\{ H_{N}(\sigma) \in A
: X_{\sigma} = 1 \}| \end{equation} converges as $N \to \infty $ to
$(\mu(A))^\ell,$ i.e. to the $\ell^{\textrm{th}}$ factorial moment
of a Poisson random variable with parameter $\mu(A)$ and the
$\ell^{\textrm{th}}$ factorial moment of $Z$ is defined as
$\mathbb{E}(Z)_\ell = \mathbb{E}Z(Z-1)\dots(Z-\ell+1).$ Then the
sequence of point processes $\mathcal{P}_N$ converges weakly to a
Poisson point process with intensity measure $\mu.$ For completeness
we proved this fact in Lemma~\ref{(a):gen_lem}.}

It is then a natural question to know if our sparseness assumption
is optimal. When the random cloud is denser can this universality
survive? We show that our sparseness condition is indeed optimal for
the number partitioning problem and for the SK model, and that the
universality does break down.

\begin{thm}\label{main_thm_2}

[Breakdown of Universality for the number partitioning problem]

\noindent (i) Let $\nu(r)=r$. Suppose that $\limsup
\frac{M(N)}{\sqrt{N}} < \infty.$ Then $\mathrm{P}$-almost surely,
the distribution of the point process $\mathcal{P}_N$ converges to
the distribution of a Poisson point process if and only if
$M=o(\sqrt{N}).$

\textrm{[Breakdown of Universality for the Sherrington-Kirkpatrick
model]}

\noindent (ii) Let $\nu(r)=r^2$. Suppose that $\limsup
\frac{M(N)}{N} < \frac{1}{8 \log 2}.$ Then $\mathrm{P}$-almost
surely, the distribution of the point process $\mathcal{P}_N$
converges to the distribution of a Poisson point process if and only
if $M=o(N).$
\end{thm}

We prove this theorem in Section \ref{(c)} by showing that the
second factorial moment of the point process does not converge to
the proper value. The case of pure $p$-spin models, with $p \ge 3,$
or more generally the case where $\nu'(0)=\nu''(0)=0$ differs
strongly (see Theorem \ref{(c):thm3}). The asymptotic behavior of
the first three moments is compatible with a Poissonian convergence.
Proving or disproving Poissonian convergence (or REM universality)
in this case is still open, as it is for the local REM conjecture.

The paper is organized as follows. In Section \ref{(ce)} we
establish important combinatorial estimates about maximal overlaps
of $\ell$-tuples of points on the random cloud. We give a particular
care to the case of pairs ($\ell=2$) and triples ($\ell=3$) which
are important for the breakdown of universality results. In Section
\ref{(a)} we establish the universality in the Gaussian case
(Theorem \ref{main_thm}). We then prove, in Section \ref{(c)}, the
breakdown of universality given in Theorem \ref{main_thm_2}. Finally
we extend the former results to a wide non-Gaussian setting in
Section \ref{NG}.

\bigskip

\section{Combinatorial estimates}\label{(ce)}

In this section we fix an integer $\ell \ge 1$ and study the maximal
overlap
\begin{equation}\label{(ce):R:max}
R_{\max}(\sigma^1,\dots,\sigma^\ell) = \max \limits_{1 \leq i < j
\leq \ell} |R(\sigma^i,\sigma^j)|.
\end{equation}
For fixed $N \in \mathbb{N}$ and $R \in [0,1)$ let us define the
following subsets of $S_N^{\ell}$
\begin{align}
& U_{\ell}(R) = \{(\sigma^1,\dots,\sigma^{\ell}):
R_{\max}(\sigma^1,\dots,\sigma^{\ell}) \leq R\} \label{(ce):U(R)}\\
\intertext{and} & V_{\ell}(R) = \{(\sigma^1,\dots,\sigma^{\ell}):
R_{\max}(\sigma^1,\dots,\sigma^{\ell}) = R\} \label{(ce):V(R)}.
\end{align}
More generally, let the sequence $R_N \in [0,1)$ be given
(respectively, the corresponding sequence of Hamming distances $d_N
= \frac{N}{2}(1 - R_N)$) and introduce the sequence of sets
$U_{\ell}(R_N)$ and $V_{\ell}(R_N)$ which we denote for simplicity
of notation by $U_{N,\ell}$ and $V_{N, \ell}$ respectively.
\par For a set $Y \subset S_{N}^{\ell}$ we denote by $Y^X$ its intersection
with $X^{\ell}$. In the following theorem we study the properties of
the sets
\begin{equation}
U_{N,\ell}^X = \{(\sigma^1, \dots, \sigma^\ell) \in X^\ell:
R_{\max}(\sigma^1, \dots, \sigma^\ell) \le R_N\}
\end{equation}
and
\begin{equation}
V_{N,\ell}^X = \{(\sigma^1, \dots, \sigma^\ell) \in X^\ell:
R_{\max}(\sigma^1, \dots, \sigma^\ell) = R_N\}.
\end{equation}
In order to state the main result of this section (Theorem
\ref{(ce):thm2}) we define the function
\begin{equation}
\label{(ce):J} \mathcal{J}(x) =
          \left\{
          \begin{array}{lr}
          \frac{1-x}{2} \log(1-x) +
          \frac{1+x}{2} \log(1+x) &
          \hbox{if $x \in [-1,1]$}, \\
          + \infty
          & \hbox{otherwise.}
          \end{array}
          \right.
\end{equation}
\begin{thm}\label{(ce):thm2} Let the sequence $R_N$ be such that
$N R_N^2 \to \infty.$

$(i)$ Then $\mathrm{P}$-almost surely
\begin{equation}
|U_{N,\ell}^X| = \mathrm{E}|U_{N,\ell}^X|(1+o(1)).
\end{equation}

(ii) If $R_N=o(1)$ and $M(N) \geq \log N$ then there exists
$\alpha\in[0,1)$ and $C>0,$ depending only on $\ell,$ such that
$\mathrm{P}$-almost surely
\begin{equation}
|V_{N,\ell}^X| \le C \, e^{\alpha N \mathcal{J}(R_N)}
\mathrm{E}|V_{N,\ell}^X|.
\end{equation}
\end{thm}

\begin{proof}

The proof is based on standard inequalities for i.i.d. random
variables that result from exponential Chebychev inequality. We
formulate them without proof: let $(X_i)_{1\leq i \leq n}$ be i.i.d.
Bernoulli rv's with $\mathbb{P}(X_i = 1)= 1 - \mathbb{P}(X_i = 0) =
p$ and let $Z = \sum \limits_{i=1}^n X_i .$ Then, for $t>0$,
\begin{align}
\label{(ce):(aux):ldev:1} & \mathbb{P}(Z - \mathbb{E}Z \geq t
\mathbb{E}Z) \leq e^{-n\left(
p(1+t)\log(1+t)+(1-p(1+t)) \log \frac{1-p(1+t)}{1-p} \right)}, \\
\label{(ce):(aux):ldev:2} & \mathbb{P}(Z - \mathbb{E}Z \leq -t
\mathbb{E}Z) \leq e^{-n\left( p(1-t)\log(1-t)+(1-p(1-t)) \log
\frac{1-p(1-t)}{1-p} \right)}.
\end{align}
If $p=p(n) \to 0$ and $t=t(n) \to 0$ as $n \to \infty,$ the above
inequalities imply that, for large enough $n$,
\begin{equation}\label{(ce):(aux):ldev:3}
\mathbb{P}\bigl( | Z - \mathbb{E}Z | \geq t \mathbb{E}Z \bigr) \leq
2 e^{-npt^2/4},
\end{equation}
whereas if $p(n) \to 0, t(n) \to \infty$, and $p(n)t(n) \to 0$ we
get from \eqref{(ce):(aux):ldev:1} that for large enough~$n$
\begin{align}\label{(ce):(aux):ldev:4}
\mathbb{P}(Z - \mathbb{E}Z \geq t \mathbb{E}Z ) \leq
e^{-\frac{1}{2}np(1+t) \log (1+t)}.
\end{align}

The proof of Theorem \ref{(ce):thm2} relies on the following
elementary lemma that again we state without proof.
\begin{lem}\label{(ce):thm1}
(i) For any sequence $R_N \in [0,1)$
\begin{equation}\label{(ce):thm1:TR}
|U_{N,\ell}| \geq 2^{N \ell}\left( 1 - 2 \binom{\ell}{2} e^{-
\frac{1}{8} N R_N^2 }\right).
\end{equation}

\noindent (ii) Suppose $R_N$ satisfies $N R_{N}^2 \to \infty$ and
$R_N=o(1).$ Then for some $C>0$ depending only on $\ell,$
\begin{equation}\label{(ce):thm1:2}
|V_{N,\ell}| = 2^{N\ell} \frac{C}{\sqrt{N}}
e^{-N\mathcal{J}(R_{N})}(1+o(1)).
\end{equation}
\end{lem}

\noindent As an elementary consequence of part $(i)$ of Lemma
\ref{(ce):thm1} one can prove that:
\begin{cor}\label{(ce):cor1}
$\mathrm{P}$-a.s. $\max \limits_{\sigma, \sigma' \in X}
|R(\sigma,\sigma')|  \leq \delta_N, $ where $\delta_N \equiv
4\sqrt{\frac{M(N) \log 2}{N} + \frac{\log N}{N}}.$
\end{cor}

The proof of part $(i)$ of Theorem \ref{(ce):thm2} then proceeds as
follows. Let us first express the size of the random cloud $|X|$ as
a sum of i.i.d. random variables
\begin{equation}
|X| = \sum \limits_{\sigma \in S_N} {\mathbf 1}_{X_{\sigma}=1}.
\end{equation}
Using \eqref{(ce):(aux):ldev:3} and the assumption that $\log N =
o(2^M)$, we see that $\mathrm{P}$-almost surely $|X|$ is given by
its expected value, i.e. $|X|= 2^M (1+o(1)).$ Therefore $|X^{\ell}|
= 2^{M\ell}(1+o(1)).$

Since $U_{N,\ell}^X \subset X^{\ell}$ and
$\mathrm{E}|U_{N,\ell}^X|=p_N^\ell |U_{N,\ell}| =
2^{M\ell}(1+o(1)),$ then proving part $(i)$ of Theorem
\ref{(ce):thm2} is equivalent to proving that the set $U_{N,\ell}^X$
coincides, up to an error of magnitude $o(2^{M\ell}),$ with the set
$X^{\ell}.$ Let us rewrite $U_{N,\ell}^X$ as
\begin{equation}\label{(ce):decomp:S}
U_{N,\ell}^X = \bigcap \limits_{1 \leq j<j' \leq \ell }
(U_{N,\ell}^X)_{jj'},
\end{equation}
where we defined
\begin{align}
(U_{N}^X)_{jj'} = \Bigl\{( \sigma^{1},\dots,\sigma^{\ell}) \in
X^{\ell}: \ |R& (\sigma^{j},\sigma^{j'})| \leq R_{N} \, \Bigr\}.
\end{align}
If we prove that every set $(U_{N,\ell}^X)_{jj'}$ coincides, up to
an error of order $o(2^{M\ell} ),$  with the set $X^{\ell},$ then
the representation \eqref{(ce):decomp:S} implies part $(i)$ of the
theorem.  We therefore concentrate on proving that
\begin{equation}
|(U_{N,\ell}^X)_{jj'}|=\mathrm{E}|(U_{N,\ell})^X_{jj'}|(1+o(1)),
\quad \mathrm{P}\mbox{-a.s.}
\end{equation}

\noindent Without loss of generality we can consider the case of
$j=1$ and $j'=2.$ By definition of $(U_{N,\ell}^X)_{jj'}$ we get
\begin{align}\label{(ce):S_jk}
|(U_{N,\ell}^X)_{12}| & = \Biggl(\sum \limits_{
\substack{\sigma^1,\sigma^2 \in S_N:\\|R(\sigma^{1},\sigma^{2})|
\leq R_N }} {\mathbf 1}_{X_{\sigma^1}=1}{\mathbf 1}_{X_{\sigma^2}=1}
\Biggr) \Biggl(\sum \limits_{\sigma \in S_N}
{\mathbf 1}_{X_{\sigma}=1} \Biggl)^{\ell-2} \nonumber \\
& = \Biggl(\sum \limits_{\sigma^1 \in S_N} {\mathbf
1}_{X_{\sigma^1}=1} \sum \limits_{\sigma^2
:|R(\sigma^{1},\sigma^{2})| \leq R_N}{\mathbf
1}_{X_{\sigma^2}=1}\Biggr) \Biggl(\sum \limits_{\sigma \in S_N}
{\mathbf 1}_{X_{\sigma}=1} \Biggr)^{\ell-2}.
\end{align}

\noindent  As we already noted, the sum in the second factor of
\eqref{(ce):S_jk} concentrates on its expected value and is equal to
$2^{M(\ell-2)} (1+o(1)).$ Let us thus turn to the first factor in
\eqref{(ce):S_jk}.

Introduce the set $(U_{N,\ell}^X)_{\sigma^1} = \{ \sigma^2 \in X
:|R(\sigma^{1},\sigma^{2})| \leq R_N\}.$ Then
\begin{equation}\label{(ce):thm2:ii:TRkX:2}
|(U_{N,\ell}^X)_{\sigma^1}| = \sum \limits_{\sigma^2:
|R(\sigma^{1},\sigma^{2})| \leq R_N}{\mathbf 1}_{X_{\sigma^2}=1}
\end{equation}
The summands in this sum are i.i.d. and it follows from part $(i)$
of Lemma \ref{(ce):thm1} that their number is at least $2^{N}(1 - 2
e^{- \frac{1}{8}N R_N^2}) = 2^{N}(1+o(1)).$ Applying
\eqref{(ce):(aux):ldev:3} together with the assumption that $\log N
= o(2^M)$ we obtain from Borel-Cantelli Lemma that
$\mathrm{P}\mbox{-a.s.}$, for any $\sigma^1\in S_N$,
\begin{equation}\label{(ce):averaging:S}
|(U_{N,\ell}^X)_{\sigma^1}| =
(1+o(1))\mathrm{E}|(U_{N,\ell}^X)_{\sigma^1}|.
\end{equation}

\noindent From \eqref{(ce):S_jk}, \eqref{(ce):thm2:ii:TRkX:2} and
\eqref{(ce):averaging:S} we immediately conclude that
\begin{equation}
|(U_{N,\ell}^X)_{12}|=(1+o(1))\mathrm{E}|(U_{N,\ell}^X)_{12}|, \quad
\mathrm{P}\mbox{-a.s.}
\end{equation}
This finishes the proof of part $(i)$ of Theorem \ref{(ce):thm2}.

\bigskip

The proof of part $(ii)$ is quite similar to the proof of part
$(i)$. By definition of $V_{N,\ell}^X$ we get
\begin{equation}\label{(ce):thm2:ii:inclusion}
V_{N,\ell}^X \subseteq \bigcup \limits_{1 \le j < j' \leq \ell}
(V_{N,\ell}^X)_{jj'},
\end{equation}
where
\begin{align}
(V_{N,\ell}^X)_{jj'} = \Bigl\{(\sigma^{1},\dots,\sigma^{\ell}) \in
V_{N,\ell}^X: \ |R&(\sigma^{j},\sigma^{{j'}})| = R_N \Bigr\}.
\end{align}

\noindent We claim that it suffices to prove that
$\mathrm{P}$-almost surely
\begin{equation} \label{(ce):thm2:ii:TRk:jj}
|(V_{N,\ell}^X)_{jj'}| \leq e^{ \alpha N \mathcal{J}(R_N)}
\mathrm{E} |(V_{N,\ell}^X)_{jj'}|.
\end{equation}
Indeed, from \eqref{(ce):thm2:ii:inclusion} and from the above
inequality we obtain that
\begin{align}\label{(ce):thm2:ii:TRkX:bnd}
|V_{N,\ell}^X| \leq\sum \limits_{1 \leq j < j' \leq \ell}
|(V_{N,\ell}^X)_{jj'}| \leq \sum \limits_{1 \leq j < j' \leq \ell}
e^{ \alpha N \mathcal{J}(R_N)} \mathrm{E}|(V_{N,\ell}^X)_{jj'}|.
\end{align}
Using part $(i)$ of Lemma \ref{(ce):thm1} it is easy to establish
that for all $1 \leq j < j' \leq \ell$
\begin{equation}
|V_{N,\ell}| = |(V_{N,\ell})_{jj'}| \bigl(1+o(1)\bigr),
\end{equation}
and therefore
\begin{equation}\label{(ce):thm2:ii:TRkX:bnd1}
\mathrm{E}|V_{N,\ell}^X| =
\mathrm{E}|(V_{N,\ell}^X)_{jj'}| (1+o(1)).
\end{equation}
Since \eqref{(ce):thm2:ii:TRkX:bnd} and
\eqref{(ce):thm2:ii:TRkX:bnd1} imply the result we concentrate on
the proof of \eqref{(ce):thm2:ii:TRk:jj}.

\par Without loss of generality we can take $j=1$ and $j'=2.$ Then,
by definition of $(V_{N,\ell}^X)_{jj'}$, we get
\begin{align}\label{(ce):(aux)_T12}
|(V_{N,\ell}^X)_{12}| & = \Biggl( \sum
\limits_{\substack{\sigma^1,\sigma^2 \in
S_N:\\|R(\sigma^{1},\sigma^{2})| = R_N }} {\mathbf
1}_{X_{\sigma^1}=1}{\mathbf 1}_{X_{\sigma^2}=1} \Biggr) \Biggl(\sum
\limits_{\sigma \in S_N} {\mathbf 1}_{X_{\sigma}=1}
\Biggr)^{\ell-2} \nonumber \\
& = \Biggl(\sum \limits_{\sigma^1 \in S_N} {\mathbf
1}_{X_{\sigma^1}=1} \sum \limits_{\sigma^2
:|R(\sigma^{1},\sigma^{2})|=R_N}{\mathbf 1}_{X_{\sigma^2}=1}\Biggr)
\Biggl(\sum \limits_{\sigma \in S_N} {\mathbf 1}_{X_{\sigma}=1}
\Biggr)^{\ell-2}.
\end{align}

\noindent As in the proof of part $(i)$ we see that the second part
of \eqref{(ce):(aux)_T12} concentrates on its expected value and
equals to $2^{M(\ell-2)} (1+o(1)).$ We are thus left to treat the
first part. Introducing the set
\begin{equation}\label{(ce):thm2:ii:VNX:1}
(V_{N,\ell}^X)_{\sigma^1} = \{ \sigma^2 \in X
:|R(\sigma^{1},\sigma^{2})|=R_N \}
\end{equation}
it is clear that
\begin{equation}\label{(ce):thm2:ii:VNX:2}
|(V_{N,\ell}^X)_{\sigma^1}| = \sum \limits_{\sigma^2
:|R(\sigma^{1},\sigma^{2})|=R_N}{\mathbf 1}_{X_{\sigma^2}=1}
\end{equation}

\par There are $2 \binom{N}{d_N}$ i.i.d. terms in the above sum.
Applying \eqref{(ce):(aux):ldev:4} with $t+1 = e^{\alpha N
\mathcal{J}(R_N)}$ where $\alpha \in [0,1)$ will be chosen later, we
obtain
\begin{equation}\label{(ce):eq25}
\mathrm{P}(|(V_{N,\ell}^X)_{\sigma^1}| \geq e^{\alpha N
\mathcal{J}(R_N)}\mathrm{E}|(V_{N,\ell}^X)_{\sigma^1}|) \leq
e^{-\frac{1}{2} \mathrm{E}|(V_{N,\ell}^X)_{\sigma^1}| (1+t) \log
(1+t)}.
\end{equation}

\noindent From Stirling's formula we see that for some $C>0$,
\begin{equation}
\mathrm{E} |(V_{N,\ell}^X)_{\sigma^1}| =  2 p_N \binom{N}{d_N} =
2^{M} \frac{C}{\sqrt{N}} e^{-N\mathcal{J}(R_N)}(1+o(1)),
\end{equation}
and thus the exponent in \eqref{(ce):eq25} is
\begin{align}
\frac{1}{2} \mathrm{E}|(V_{N,\ell}^X)_{\sigma^1}| & (1+t) \log (1+t)
= \frac{C}{2 \sqrt{N}} \, 2^{M} e^{-(1-\alpha) N\mathcal{J}(R_N)} \,
\alpha N \mathcal{J}(R_N)(1 + o(1)) \nonumber \\
& = \frac{\alpha C}{2}e^{ M \log 2 -(1-\alpha) N \mathcal{J}(R_N) -
\frac{1}{2}\log N } N \mathcal{J}(R_N)(1 + o(1)).
\end{align}
If for every $\sigma^1 \in S_N$ the set $(V_{N,\ell}^X)_{\sigma^1}$
is empty then there is nothing to prove. Otherwise, by Corollary
\ref{(ce):cor1} we obtain that $\mathrm{P}$-almost surely $R_N <
\delta_N$, and since $\mathcal{J}(x) = \frac{1}{2}x^2(1+O(x^2))$
near the origin we can choose $\alpha \in [0,1)$ in such a way that
\begin{align}
M \log 2 - (1-\alpha) N\mathcal{J}(R_N) - \frac{1}{2}\log N > \gamma
\log N
\end{align}
for some $\gamma > 0.$ Hence $\sum \limits_N e^{-\frac{1}{2}
\mathrm{E}|(V_{N,\ell}^X)_{\sigma^1}| (1+t) \log (1+t)} < \infty,$
and we obtain from \eqref{(ce):eq25} and Borel-Cantelli Lemma that
$\mathrm{P}$-almost surely
\begin{equation}\label{(ce):T_sigma_1_bound}
|(V_{N,\ell}^X)_{\sigma^1}| \leq e^{\alpha N \mathcal{J}(R_N)}
\mathrm{E}|(V_{N,\ell}^X)_{\sigma^1}|.
\end{equation}
It is easy to see that \eqref{(ce):T_sigma_1_bound} and
\eqref{(ce):(aux)_T12} imply \eqref{(ce):thm2:ii:TRk:jj}. This
concludes the proof of Theorem \ref{(ce):thm2}.
\end{proof}
\medskip

In part $(ii)$ of Theorem \ref{(ce):thm2} we studied the properties
of the sets $V_{N,\ell}^X \subset X^{\ell}$ for arbitrary $\ell \geq
2.$ For $\ell=2$ we can improve on Theorem~\ref{(ce):thm2}.

\begin{thm}\label{(ce):thm3} Suppose that $\limsup \frac{M}{N} < 1.$

(i) If for some $c_1 > \frac{1}{2}$
\begin{equation}\label{(ce):thm3:1}
\mathcal{J}(R_N) \leq \frac{M \log 2}{N}  - \frac{c_1 \log N}{N},
\end{equation}
then $\mathrm{P}$-almost surely
\begin{equation}
|V_{N,2}^X| = (1+o(1)) \mathrm{E}|V_{N,2}^X|.
\end{equation}

(ii) If for positive constants $c_1,c_2$
\begin{equation}\label{(ce):thm3:2}
\frac{M \log 2}{N} - \frac{c_1 \log N}{N} \leq \mathcal{J}(R_N) \leq
\frac{M \log 2}{N} + \frac{c_2 \log N}{N}
\end{equation}
then there is a constant $c$ such that $\mathrm{P}$-almost surely
\begin{equation}
|V_{N,2}^X| \leq N^c \ \mathrm{E}|V_{N,2}^X|.
\end{equation}

(iii) If for some $c_2 > \frac{3}{2}$
\begin{equation}\label{(ce):thm3:3}
\mathcal{J}(R_N) > \frac{M \log 2}{N} + \frac{c_2 \log 2}{N},
\end{equation}
then the set $V_{N,2}^X$ is $\mathrm{P}$-almost surely empty.
\end{thm}

\begin{proof}

$(i)$ By definition of $V_{N,2}^X$,
\begin{equation}\label{(ce):T:R12}
|V_{N,2}^X| = \sum \limits_{\sigma^1 \in S_N} {\mathbf
1}_{X_{\sigma^1}=1} \sum \limits_{\sigma^2 \in S_N:
|R(\sigma^1,\sigma^2)|=R_{N}} {\mathbf 1}_{X_{\sigma^2}=1}.
\end{equation}
Using \eqref{(ce):thm2:ii:VNX:1} the inner sum is
$|(V_{N,2}^X)_{\sigma^1}|.$ Since it is a sum of i.i.d. random
variables then for all $t=t(N)=o(1)$ we get from
\eqref{(ce):(aux):ldev:3} that
\begin{equation}\label{(ce):concentr:11}
\mathrm{P}\left( | |(V_{N,2}^X)_{\sigma^1}| -
\mathrm{E}|(V_{N,2}^X)_{\sigma^1}| | \geq t \, \mathrm{E}
|(V_{N,2}^X)_{\sigma^1}| \right) \leq 2 e^{-
{t^2}\mathrm{E}|(V_{N,2}^X)_{\sigma^1}|/4}.
\end{equation}

\noindent Using Stirling's approximation
\begin{equation}
\sqrt{2\pi} n^{n+\frac{1}{2}} e^{-n+\frac{1}{12n+1}}<
n!<\sqrt{2\pi}n^{n+\frac{1}{2}}e^{-n+\frac{1}{12n}}
\end{equation}
we obtain that
\begin{align}\label{(ce):thm3:i:stirling}
\mathrm{E}& |(V_{N,2}^X)_{\sigma^1}| = 2 p_N \binom{N}{d_N} =
O\Bigl( \frac{e^{M \log 2 - N \mathcal{J}(R_N)}}{\sqrt{N(1 -
R_N^2)}} \Bigr ).
\end{align}
Further, from \eqref{(ce):thm3:1}  and Corollary \ref{(ce):cor1} we
obtain that $\mathrm{E}|(V_{N,2}^X)_{\sigma^1}| \ge C N^{c_1-1/2}$
for some positive constant $C>0$. Choosing  $t = N^{-\gamma}$ for
some small enough $\gamma >0$ in \eqref{(ce):concentr:11}, we
conclude from Borel-Cantelli Lemma that
\begin{equation}\label{(ce):thm3:1:1}
|(V_{N,2}^X)_{\sigma^1}| =
\mathrm{E}|(V_{N,2}^X)_{\sigma^1}|(1+o(1)), \quad
\mathrm{P}\mbox{-a.s.}
\end{equation}
Using \eqref{(ce):T:R12}, \eqref{(ce):thm3:1:1}, and the fact that
$\mathrm{P}$-almost surely there are $2^{M}(1+o(1))$ configurations
in the random cloud $X$ proves $(i).$
\medskip

$(ii)$ The proof is similar to the proof of part $(i).$ In
particular, from the representation \eqref{(ce):T:R12} it is easy to
see that it suffices to prove that
\begin{equation}
|(V_{N,2}^X)_{\sigma^1}| \leq N^c \ \mathrm{E}
|(V_{N,2}^X)_{\sigma^1}|, \quad \mathrm{P}\mbox{-a.s.}
\end{equation}
Choosing $1+t = N^c$ we  get from \eqref{(ce):thm3:i:stirling} that
for some positive constant $C$
\begin{align}\label{(ce):thm3:ii:conctr}
\frac{1}{2} E &|(V_{N,2}^X)_{\sigma^1}| (1+t) \log (1+t) \ge C
N^{c-c_2-1/2} \log N.
\end{align}
Choosing $c$ large enough and applying \eqref{(ce):(aux):ldev:4}
together with Borel-Cantelli Lemma proves part $(ii)$.
\medskip

$(iii)$ We again use the representation \eqref{(ce):T:R12}. Clearly
it is enough to prove that for all $\sigma^1 \in X$ the set
$(V_{N,2}^X)_{\sigma^1}$ is $\mathrm{P}$-almost surely empty. By
\eqref{(ce):thm3:i:stirling} $\mathrm{E}|(V_{N,2}^X)_{\sigma^1}| \le
C N^{-c_2 - 1/2}$ for some positive constant $C.$ Thus, choosing
$1+t = \frac{1}{E|(V_{N,2}^X)_{\sigma^1}|} \to \infty$ we obtain
from~\eqref{(ce):(aux):ldev:4}
\begin{align}
\mathrm{P} & \left( |(V_{N,2}^X)_{\sigma^1}| \geq 1\right) =
\mathrm{P} \left(|(V_{N,2}^X)_{\sigma^1}| \geq
(1+t)\mathrm{E}|(V_{N,2}^X)_{\sigma^1}| \right) \leq (t+1)^{-1/2}.
\end{align}
By definition of $t$ and by condition \eqref{(ce):thm3:3} we get
\begin{align}
(t+1)^{-1/2} \leq C^{1/2} N^{-c_2/2 - 1/4}.
\end{align}
Applying Borel-Cantelli Lemma proves $(iii)$.
\end{proof}
\medskip

In order to estimate the third moment in Theorems \ref{(c):thm2},
\ref{(c):thm1}, and \ref{(c):thm3} we give a result similar to
Theorem \ref{(ce):thm3} but for $\ell=3,$ i.e. for a given sequence
of vectors $(R_{12}^N, R_{23}^N, R_{31}^N)$ we estimate the cardinal
of the set
\begin{equation}\label{(ce):W}
W_{N,3}^X = \Bigl\{(\sigma^1, \sigma^2, \sigma^3) \in X^{3}:
R(\sigma^1,\sigma^2)=R_{12}^N, R(\sigma^2,\sigma^3)=R_{23}^N,
R(\sigma^3,\sigma^1)=R_{31}^N \Bigr\}.
\end{equation}
Below we omit the explicit dependence of the sequence $(R_{12}^N,
R_{23}^N, R_{31}^N)$ on $N$ and instead of $R_{12}^N, R_{23}^N$ and
$R_{31}^N$ will write  $R_{12}, R_{23}$ and $R_{31}$ respectively.
In order to formulate the theorem we introduce the following
function on $\mathbb{R}^3:$
\begin{align}
\mathcal{J}^{(2)}(& x,y,z) = \frac{1+ x + y + z}{4} \log (1 + x + y
+
z )  + \frac{1 + x - y - z}{4} \log (1 + x - y - z ) \nonumber \\
& + \frac{1 - x + y - z}{4} \log (1 - x + y - z ) + \frac{1 - x - y
+ z}{4} \log (1 - x - y + z )
\end{align}
if $|1 + x| \ge |y + z|$ and $|1 - x | \ge |y - z|$ and
$\mathcal{J}^{(2)}(x,y,z) = + \infty$ otherwise.

\begin{thm}\label{(ce):thm6} Suppose $\limsup \frac{M}{N} < 1.$

$(i)$ If for some $c_1>\frac{1}{2}$ the sequence $R_{12}$ satisfies
\eqref{(ce):thm3:1}, and if for some $c_{1}^{(2)} > 1$
\begin{equation}\label{(ce):thm6:1}
\mathcal{J}^{(2)}(R_{12},R_{23},R_{31}) \leq \frac{M \log 2}{N} +
\mathcal{J}(R_{12}) - \frac{c_1^{(2)} \log N}{N},
\end{equation}
then $\mathrm{P}$-almost surely
\begin{equation}
|W_{N,3}^X | = (1+o(1))\mathrm{E}|W_{N,3}^X|.
\end{equation}

$(ii)$ If for positive constants $c_1^{(2)},c_2^{(2)}$,
\begin{align}\label{(ce):thm6:2}
\frac{M \log 2}{N} + \mathcal{J}(R_{12}) - \frac{c_1^{(2)} \log
N}{N} & \leq \mathcal{J}^{(2)}(R_{12},R_{23},R_{31}) \nonumber \\
& \leq \frac{M \log 2}{N} + \mathcal{J}(R_{12}) + \frac{c_2^{(2)}
\log N}{N},
\end{align}
then there is a constant $c$ such that $\mathrm{P}$-almost surely
\begin{equation}
|W_{N,3}^X| \leq N^c \mathrm{E}|W_{N,3}^X|.
\end{equation}

$(iii)$ If for some $c_2^{(2)} > \frac{3}{2}$
\begin{equation}\label{(ce):thm6:3}
\mathcal{J}^{(2)}(R_{12},R_{23},R_{31}) > \frac{M \log 2}{N}  +
\mathcal{J}(R_{12}) + \frac{c_2^{(2)} \log 2}{N},
\end{equation}
then $\mathrm{P}$-almost surely the set $W_{N,3}^X$ is empty.
\end{thm}

\begin{proof}

>From Lemma \ref{(ce):p:2:m:3:aux1} stated below it follows that, for
arbitrary configurations $\sigma^1, \sigma^2,\sigma^3\in S_N$, the
function $\mathcal{J}^{(2)}(R_{12},R_{23},R_{31})$ is well defined.
This lemma, whose proof we will omit, is a direct consequence of the
fact that the Hamming distance on $S_N$ satisfies triangle
inequality.
\begin{lem}\label{(ce):p:2:m:3:aux1}
For arbitrary configurations  $\, \sigma^1, \sigma^2$ and $\sigma^3
\in S_N$ we have $1 + R_{12} \ge |R_{23} + R_{31}|$ and $1 - R_{12}
\ge |R_{23} - R_{31}|.$
\end{lem}

The proof of the Theorem \ref{(ce):thm6} is similar to that of
Theorem \ref{(ce):thm3}. We begin by writing the size of the set
$W_{N,3}^X$ as
\begin{equation}\label{(b):p:2:m:3:S:X}
|W_{N,3}^X| = \sum \limits_{ \sigma^1 \in S_N} {\bf
1}_{X_{\sigma^1}=1} \sum \limits_{\substack{\sigma^2 \in
S_N:\\R(\sigma^1, \sigma^2)=R_{12}}} {\bf 1}_{X_{\sigma^2}=1} \sum
\limits_{\substack{\sigma^3 \in
S_N:\\R(\sigma^2,\sigma^3)=R_{23}\\R(\sigma^1,\sigma^3)=R_{31}}}
{\bf 1}_{X_{\sigma^3}=1}
\end{equation}
Let us first estimate the number of terms in the last sum. This
means that, given $\sigma^1, \sigma^2 \in S_N$ with overlap
$R(\sigma^1, \sigma^2) = R_{12}$, we have to calculate the number of
configurations $\sigma^3$ with $R(\sigma^2,\sigma^3)=R_{23}$ and
$R(\sigma^3,\sigma^1)=R_{31}.$ Without loss of generality we can
assume that  all the spins of $\sigma^1$ are equal to $1.$ Further,
let $C(\sigma^1,\sigma^2,\sigma^3)$ be a $3 \times N$ matrix with
rows $\sigma^1,\sigma^2,\sigma^3.$ For a column vector $\boldsymbol
\delta \in \{-1,1\}^3$ we let $n_{\boldsymbol \delta}$ be the number
of columns of the matrix $C$ that are equal to $\boldsymbol \delta,$
i.e.
\begin{equation}
n_{\boldsymbol \delta} = \bigl|\{j \leq N: \, (\sigma^1_j,
\sigma^2_j, \sigma^3_j) = {\boldsymbol \delta}\}\bigr|.
\end{equation}
Then the overlaps can be written in terms of $n_{\boldsymbol
\delta},$ namely
\begin{equation}
\left\{
\begin{array}{lr}
n_{(1,1,1)} + n_{(1,1,-1)} - n_{(1,-1,1)} - n_{(1,-1,-1)} = NR_{12},\\
n_{(1,1,1)} - n_{(1,1,-1)} - n_{(1,-1,1)} + n_{(1,-1,-1)} = NR_{23},\\
n_{(1,1,1)} - n_{(1,1,-1)} + n_{(1,-1,1)} - n_{(1,-1,-1)} = NR_{31},\\
n_{(1,1,1)} + n_{(1,1,-1)} + n_{(1,-1,1)} + n_{(1,-1,-1)} = N.
\end{array}
\right.
\end{equation}
Solving this system of linear equations we find
\begin{equation}
\left\{
\begin{array}{ll}
n_{(1,1,1)}  & = \frac{1}{4}N(1 + R_{12} + R_{23} + R_{31}), \\
n_{(1,1,-1)} & = \frac{1}{4}N(1 + R_{12} - R_{23} - R_{31}), \\
n_{(1,-1,1)} & = \frac{1}{4}N(1 - R_{12} - R_{23} + R_{31}), \\
n_{(1,-1,-1)} & = \frac{1}{4}N(1 - R_{12} + R_{23} - R_{31}). \\
\end{array}
\right.
\end{equation}
We notice that specifying the configuration $\sigma^3$ is equivalent
to specifying the numbers $n_{(1,1,1)},n_{(1,1,-1)},n_{(1,-1,1)}$,
and $n_{1,-1,-1}.$ Therefore the number of configurations $\sigma^3
\in S_N$ with overlaps $R(\sigma^2,\sigma^3)=R_{23}$, and
$R(\sigma^3,\sigma^1)=R_{31}$ is
\begin{align}\label{(ce):thm6:4}
& \binom{n_{(1,1,1)} + n_{(1,1,-1)}}{n_{(1,1,1)}}
\binom{n_{(1,-1,1)}
+ n_{(1,-1,-1)}}{n_{(1,-1,1)}} \nonumber \\
& \qquad = \frac{\bigl(n_{(1,1,1)} +
n_{(1,1,-1)}\bigr)!}{n_{(1,1,1)}! \, n_{(1,1,-1)}!}
\frac{\bigl(n_{(1,-1,1)} + n_{(1,-1,-1)}\bigr)!}{n_{(1,-1,1)}! \,
n_{(1,-1,-1)}!}.
\end{align}
Applying Stirling's approximation to \eqref{(ce):thm6:4} one obtains
that the number of terms in the last  summation in
\eqref{(b):p:2:m:3:S:X} is of order
$$
\frac{2^N e^{N \mathcal{J}(R_{12}) - N
\mathcal{J}^{(2)}(R_{12},R_{23},R_{31})} \sqrt{1 - R_{12}^2}}{ N
\sqrt{P(R_{12},R_{23},R_{31})}},
$$ where
$P(x,y,z)=(1+ x + y + z)(1+x - y - z)(1 - x + y - z)(1 - x - y +
z).$ The rest of the proof essentially is a rerun of the proof of
Theorem \ref{(ce):thm3}. We skip the details.
\end{proof}

\medskip

{\section{Proof of Theorem \ref{main_thm}}\label{(a)}}

As in \cite{BK} and \cite{BCMN1},\cite{BCMN2}, the proof of the
Poisson convergence is based on the analysis of factorial moments of
the point processes $\mathcal{P}_N$ defined in \eqref{1}.

In general, let $\xi_N$ be a sequence of point processes defined on
a common probability space $(\Omega, \mathcal{F}, \mathbb{Q})$ and
let $\xi$ be a Poisson point process with intensity measure $\mu.$
Define the $\ell^{\textrm{th}}$ factorial moment
$\mathbb{E}(Z)_\ell$ of the random variable $Z$ to be
$\mathbb{E}Z(Z-1)\dots(Z-\ell+1).$ The following is a classical
lemma that is a direct consequence of Theorem 4.7 in
\cite{Kallenberg}.
\medskip

\begin{lem}\label{(a):gen_lem}

If for every $\ell \ge 1$ and every Borel set $A$
\begin{equation}\label{(a):PPP_conv}
\lim \limits_{N \to \infty} \mathbb{E_Q}(\xi_N(A))_\ell =
(\mu(A))^\ell,
\end{equation}
then the distribution of $(\xi_N)_{N\geq 1}$ converges weakly to the
distribution of $\xi.$
\end{lem}

Applying Lemma \ref{(a):PPP_conv} to the sequence of point processes
$\mathcal{P}_N = \sum \limits_{\sigma \in X} \delta_{H'_N(\sigma)}$
the following result proves Theorem \ref{main_thm}.

\begin{thm}\label{lem2}
Under the assumptions of Theorem \ref{main_thm}, for every $\ell \in
\mathbb{N}$ and every bounded Borel set $A$
\begin{equation}\label{(a):fact:mom}
\lim \limits_{N \to \infty} \mathbb{E}(\mathcal{P}_N(A))_\ell =
(\mu(A))^\ell, \quad \mathrm{P}\mbox{-a.s.},
\end{equation}
where $\mu$ is defined in \eqref{(intro):intensity}.
\end{thm}

\begin{proof}

\par We start with the computation of the first moment of $\mathcal{P}_N(A):$
\begin{equation}\label{(a):g:first:mom}
\mathbb{E}\, \mathcal{P}_N(A) = \sum \limits_{\sigma \in X}
\mathbb{P}(H'_N(\sigma) \in A).
\end{equation}
As we saw in the proof of Theorem \ref{(ce):thm2} the size of the
random cloud $|X|$ is $\mathrm{P}$-almost surely $2^M(1+o(1)).$
Since $H'_N(\sigma), \sigma \in S_N$, are identically distributed
normal random variables with mean $-{a_N}/{b_N}$ and variance
${1}/{b_N}$, the sum in \eqref{(a):g:first:mom} can be written as
\begin{equation}\label{(a):first:momnt}
2^M \frac{e^{-a_N^2/{2}}}{\sqrt{2 \pi}} b_N \int \limits_{A} e^{-x^2
b_N^2/2 - a_N b_N x} dx \,(1+o(1)), \quad P\mbox{-a.s.}
\end{equation}
By the dominated convergence theorem and the definition of $a_N$ it
follows from \eqref{(a):first:momnt} that the limit of the first
moment is $\mu(A).$

\par To calculate factorial moments of higher order we follow
\cite{BCMN1} and rewrite the $\ell^{\mathrm{th}}$ factorial moment
of $\mathcal{P}_N(A)$ as
\begin{equation}\label{(a):gen_mom}
\mathbb{E} ( \mathcal{P}_N(A) )_{\ell} = \sum
\limits_{\sigma^1,\dots,\sigma^{\ell}} \mathbb{P}\bigl(
H'_N(\sigma^1) \in A, \dots, H'_N(\sigma^\ell) \in A\bigr)
\end{equation}
where the sum runs over all ordered sequences of distinct
configurations $(\sigma^1,\dots,\sigma^{\ell})\in X^{\ell}.$ To
analyze it we decompose the set $S_N^{\ell}$ into three
non-intersecting subsets
\begin{equation}\label{(a):decomposition}
S_N^{\ell}= U_\ell (R_N) \cup \bigl( U_\ell (\delta_N) \backslash
U_\ell (R_N)\bigr) \cup (U_\ell (\delta_N))^c,
\end{equation}
where $U_\ell$ is defined in \eqref{(ce):U(R)}, $\delta_N$ is
defined in Corollary \ref{(ce):cor1} and the sequence $R_N$ is
chosen such that:
\begin{equation}\label{(a):g:c_N}
R_N \to 0, M \nu (R_N) \to 0,  N R_N^2 \to \infty.
\end{equation}
This is possible since we have assumed that $M=o(\sqrt{N})$ if
$\nu'(0) \neq 0$ and $M=o(N)$ if $\nu'(0) = 0$. We recall here
that the function $\nu$ is defined in
\eqref{(intro):cov:str}-\eqref{(intro):cov:str'}.

\par Having specified $R_N$, let us analyze the contribution to the sum
\eqref{(a):gen_mom} coming from the intersection of $X^{\ell}$ with
the sets $U_\ell(R_N), U_\ell (\delta_N)\backslash U_\ell (R_N)$,
and $(U_\ell (\delta_N))^c.$ Firstly, by Corollary \ref{(ce):cor1}
the intersection of the set $(U_\ell (\delta_N))^c$ with $X^{\ell}$
is $\mathrm{P}$-a.s. empty and therefore its contribution to the sum
\eqref{(a):gen_mom} is zero. Next, let us show that
$\mathrm{P}$-a.s.
\begin{align}\label{(a):contr:S}
& \lim \limits_{N \to \infty} \sum
\limits_{\sigma^{1},\dots,\sigma^{\ell} } \mathbb{P} \left(
H'_N(\sigma^{1}) \in A,\dots, H'_N(\sigma^{\ell}) \in A \right) =
({\mu(A)})^\ell,
\end{align}
where the sum is over all the sequences
$(\sigma^1,\dots,\sigma^\ell) \in U_\ell (R_N) \cap X^{\ell}.$

\par For every $\ell \in \mathbb{N}$ let $B(\sigma^1,\dots,\sigma^\ell)$
denote the covariance matrix of  random variables
$H_N(\sigma^1),\dots,H_N(\sigma^\ell).$ By \eqref{(intro):cov:str}
its elements, $b_{ij},$ are given by
\begin{equation}
b_{ij} = \nu(R_{ij})
\end{equation}
where we wrote $R(\sigma^i,\sigma^j)=R_{ij}.$ Since $R_N = o(1)$ the
matrix $B(\sigma^1,\dots,\sigma^\ell)$ is non-degenerate. We
therefore get for $(\sigma^{1}, \dots,\sigma^{\ell}) \in U_\ell
(R_N) \cap X^{\ell}$ that
\begin{align}\label{(a):explicit:formula}
\mathbb{P} \bigl( H'_{N}&(\sigma^{1}) \in
A,\dots,H'_{N}(\sigma^{\ell} ) \in A \bigr) =
\frac{b_N^\ell}{(2\pi)^{\ell/2}\sqrt{\det B}} \nonumber \\ & \times
\int \limits_{A} \dots \int \limits_{A} e^{- ( \vec{x}, B^{-1}
\vec{x} ) b_N^2 /2 - a_N b_N ( \vec{x} , B^{-1} \vec{1} ) - a_{N}^2
( \vec{1}, B^{-1} \vec{1} ) / 2 } d \vec{x},
\end{align}
where $B = B(\sigma^{1},\dots,\sigma^{\ell}), \vec{x} =
(x_{1},\dots,x_{\ell})^{t}$ and $\vec{1}=(1,\dots,1)^{t}.$ From the
definition of $a_N$ we further get from \eqref{(a):explicit:formula}
that
\begin{align}\label{(a):common_term}
\mathbb{P} \bigl( H'_{N}& (\sigma^{1}) \in A,\dots,H'_{N}
(\sigma^{\ell}) \in A \bigr) = \frac{1}{2^{M\ell}} \frac{e^{(\ell -
( \vec{1}, B^{-1} \vec{1} ))(
M\log 2 - 1/2 \log M)}} {\sqrt{ \det B} } \nonumber \\
& \times \frac{1}{ (\sqrt{\pi})^{\ell} }\int \limits_{A} \dots \int
\limits_{A} e^{-  ( \vec{x}, B^{-1} \vec{x} ) b_N^2 / 2  - a_N b_N (
\vec{x}, B^{-1}\vec{1}) } d \vec{x}.
\end{align}

\par Since the matrix $B^{-1}$ is positive definite
and since $a_N b_N \to \sqrt{2 \log 2}$ we conclude from the
dominated convergence theorem that for all bounded Borel sets $A,$
uniformly in $(\sigma^{1}, \dots,\sigma^{\ell}) \in U_\ell (R_N)
\cap X^{\ell},$
\begin{align}\label{(a):int_part}
\frac{1}{( \sqrt{\pi})^{\ell}} & \int \limits_{A} \dots \int
\limits_{A} e^{ -(\vec{x}, B^{-1} \vec{x})b^2_N/2 - a_N b_N
(\vec{x}, B^{-1} \vec{1} ) } d \vec{x} \nonumber \\ & \to \frac{1}{(
\sqrt{\pi})^{\ell}} \int \limits_{A} \dots \int \limits_{A} e^{ -
\sqrt{2 \log 2} (\vec{x},\vec{1} ) } d \vec{x} = (\mu(A))^\ell.
\end{align}

\par To evaluate $e^{(\ell - (\vec{1}, B^{-1} \vec{1} ))( M\log 2 -
1/2 \log M)} $ in \eqref{(a):common_term} we look at
$(\ell-(\vec{1},B^{-1}\vec{1}) )\det B$ as a multivariate function
of $(b_{ij})_{1\leq i < j \leq \ell}.$ It is a polynomial of degree
$\ell$ with coefficients depending only on $\ell$ and without
constant term. It implies that $|\ell-(\vec{1},B^{-1} \vec{1})| =
O(\nu(R_{N}))$ and therefore, by \eqref{(a):g:c_N},
\begin{equation}\label{(a):exp_part}
e^{( \ell-(\vec{1},B^{-1} \vec{1}) )(M \log 2 - 1/2 \log M)} = 1 +
o(1).
\end{equation}

\noindent Combining \eqref{(a):common_term}, \eqref{(a):int_part},
and \eqref{(a):exp_part} we may rewrite the sum
\eqref{(a):contr:S}~as
\begin{align}\label{(a):sum_S}
\sum \limits_{\sigma^{1},\dots,\sigma^{\ell}} \frac{1}{2^{M\ell}}
\left( \mu(A)\right)^\ell (1+o(1))= \frac{|U_\ell (R_N) \cap
X^\ell|}{2^{M\ell}} \left( \mu(A)\right)^\ell (1+o(1)).
\end{align}

Now it follows from Theorem \ref{(ce):thm2} $(i)$ and Lemma
\ref{(ce):thm1} $(i)$ that $|U_\ell (R_N) \cap X^\ell|$ concentrates
around its expected value, namely $2^{M \ell}(1+o(1)),$ so that
\eqref{(a):sum_S} implies~\eqref{(a):contr:S}.

Next, let us establish that the contribution from the second set in
\eqref{(a):decomposition} is negligible, i.e. let us prove that
\begin{align}\label{(a):contr:T2}
& \lim \limits_{N \to \infty} \sum
\limits_{\sigma^{1},\dots,\sigma^{\ell}}\mathbb{P} \left(
H'_N(\sigma^{1}) \in A,\dots, H'_N(\sigma^{\ell}) \in A \right) = 0,
\end{align}
where the sum runs over all the sequences
$(\sigma^1,\dots,\sigma^\ell) \in \bigl( U_\ell (\delta_N)
\backslash U_\ell (R_N)\bigr) \cap X^{\ell}.$ To do this we first
bound the righthand side of \eqref{(a):common_term}. By definition
\eqref{(ce):R:max} and Corollary~\ref{(ce):cor1},
$R_{\max}(\sigma^1,\dots,\sigma^\ell) \leq \delta_N = o(1).$
Therefore following the same reasoning as above we obtain from
\eqref{(a):common_term} that for some constant $C>0$ and for all
$(\sigma^1,\dots,\sigma^\ell) \in \bigl( U_\ell(\delta_N) \backslash
U_\ell (R_N)\bigr) \cap X^{\ell}$
\begin{equation}\label{(a):g:gen_term_bound}
\mathbb{P} \bigl( H'_{N}(\sigma^{1}) \in A,\dots,H'_{N}
(\sigma^{\ell}) \in A \bigr) \leq \frac{e^{CM \nu( R_{\max}
)}}{2^{M\ell}} (\mu(A))^{\ell}(1+o(1)).
\end{equation}
For fixed $N$ the overlap takes only a discrete set of values
\begin{equation}
K_N = \Bigl\{1 - \frac{2k}{N}: k=0,1,\dots,N \Bigr\}.
\end{equation}
We represent the set $U_\ell(\delta_N) \backslash U_\ell (R_N)$ as a
union of sets $V_{\ell}(R_{N,k}),$ where we denoted $R_{N,k} = 1 -
\frac{2k}{N} \in K_N \cap (R_N, \delta_N ).$ Let us fix $k$ and
bound the contribution from the set $V_{\ell}(R_{N,k}) \cap X^\ell,$
i.e.
\begin{align}\label{(a):g:T_2^k_sum}
& \sum \limits_{ \substack{(\sigma^{1},\dots,\sigma^{\ell}) \in \\
V_{\ell}(R_{N,k}) \cap X^\ell} } \mathbb{P} \left( H'_N(\sigma^{1})
\in A,\dots, H'_N(\sigma^{\ell}) \in A \right).
\end{align}
\noindent We obtain from Lemma \ref{(ce):thm1} $(ii)$ that
\begin{equation}
|V_{\ell}(R_{N,k})| = 2^{N \ell} \frac{C}{\sqrt{N}} e^{-N
\mathcal{J}(R_{N,k})}(1+o(1)).
\end{equation}

In the case $M(N) \leq \log N$ we can choose the sequence $R_N$ in
such a way that the set $U_\ell(\delta_N) \backslash U_\ell (R_N)$
is empty. Therefore we can assume without loss of generality that
$M(N) \geq \log N.$ Applying part $(ii)$ of Theorem \ref{(ce):thm2}
we further conclude that
\begin{equation}
|V_{\ell}(R_{N,k}) \cap X^\ell| \leq 2^{M \ell }\frac{C}{\sqrt{N}}
e^{-(1-\alpha)N \mathcal{J}(R_{N,k})}, \quad \mathrm{P}\mbox{-a.s.}
\end{equation}
Using \eqref{(a):g:gen_term_bound} and the last inequality we bound
the sum \eqref{(a):g:T_2^k_sum} by
\begin{align}\label{(a):g:T_2^k_bound1}
& \ \frac{C}{\sqrt{N}} e^{-(1-\alpha) N \mathcal{J}(R_{N,k})} e^{C M
\nu(R_{N,k})}.
\end{align}

\par One can easily check that $M \nu(R_{N,k}) = o(N R_{N,k}^2)$
for $R_{N,k} \in (R_N, \delta_N ).$ Together with $\mathcal{J}(x)
\geq x^2/2$ it implies that for some positive constants $C_1$ and
$C_2$ we can further bound \eqref{(a):g:T_2^k_bound1} by
\begin{align}\label{contr from T_2^k}
\frac{C_1}{\sqrt{N}} e^{-C_2 N R_{N,k}^2}.
\end{align}

\noindent As a consequence, we obtain an almost sure bound
\begin{align}
\sum \limits_{ \substack{R_{N,k} \in \\ K_N \cap (R_N, \delta_N)}} &
\sum
\limits_{\substack{(\sigma^{1},..,\sigma^{\ell}) \in \\
V_{\ell}(R_{N,k}) \cap X^\ell}} \mathbb{P} \left(
H'_N(\sigma^{1}) \in A,\dots, H'_N(\sigma^{\ell}) \in A \right) \nonumber\\
& \leq \frac{C_1}{\sqrt{N}} \sum \limits_{\substack{R_{N,k} \in \\
K_N \cap (R_N, \delta_N)}} e^{-C_2 N R_{N,k}^2}.
\end{align}

\noindent Introducing new variables $y_{N,k} = \sqrt{N} R_{N,k}$ we
rewrite the above sum as
\begin{equation}
\frac{C_1}{2} \sum \limits_{y_{N,k}} e^{-C_2 y_{N,k}^2}
\frac{2}{\sqrt{N}},
\end{equation}
where the summation is over the discrete set $\sqrt{N} K_N \, \cap
(\sqrt{N}R_N,\sqrt{N} \delta_N).$ Since $NR_N^2 \to \infty,$ then for
arbitrary $C>0$ and for large $N$ we can further bound this sum by
\begin{equation}
\frac{C_1}{2} \sum \limits_{y_{N,k} \geq C } e^{-C_2 y_{N,k}^2}
\frac{2}{\sqrt{N}}.
\end{equation}
Interpreting the last sum as a sum of areas of nonintersecting
rectangles with one side equal $e^{-C_2 y_{N,k}}$ and the other
$\frac{2}{\sqrt{N}}$ we bound it with the integral
\begin{equation}
\int \limits_{C - \frac{2}{\sqrt{N}}}^{\infty} e^{-C_2 y^2} dy
\end{equation}
Since the constant $C$ is arbitrary we get that \eqref{(a):contr:T2}
is $o(1).$ This finishes the proof of Theorem~\ref{lem2} and
therefore of Theorem \ref{main_thm}.
\end{proof}

\bigskip

\section{Proof of Theorem \ref{main_thm_2}.}\label{(c)}

In order to prove the breakdown of universality in
Theorem~\ref{main_thm_2} we use a strategy similar to that used in
\cite{BCMN2} to disprove the local REM conjecture for the number
partitioning problem and for the Sherrington-Kirkpatrick model when
the energy scales are too large. We prove that $\mathrm{P}$-a.s. for
every bounded Borel set $A$
\begin{enumerate}
  \item the limit of the first factorial moment exists and equals $\mu(A);$
  \item the second factorial moment $\mathbb{E}(\mathcal{P}_N(A))_{2}$ does not converge
  to $(\mu(A))^2;$
  \item the third moment is bounded.
\end{enumerate}
These three facts immediately imply that the sequence of random
variables $\mathcal{P}_N(A)$ does not converge weakly to a Poisson
random variable and so the sequence of point processes
$\mathcal{P}_N$ does not converge weakly to a Poisson point process.
Part $(i)$ of Theorem \ref{main_thm_2} is thus obviously implied by
the following

\begin{thm}[Breakdown of Universality for the number partitioning problem]\label{(c):thm2}
Let $\nu(r)=r$. For every bounded Borel set A
\begin{gather}
\lim \mathbb{E} (\mathcal{P}_N(A))_1 = \mu(A), \quad
\mathrm{P}\mbox{-a.s.}
\end{gather}
Moreover, if $\limsup \frac{M(N)}{\sqrt{N}} = \varepsilon < \infty$
then $\mathrm{P}$-a.s.

(i) $\limsup \mathbb{E} (\mathcal{P}_N(A))_2 = e^{2 \varepsilon^2
\log^2 2}(\mu(A))^2,$

(ii) $\limsup \mathbb{E} (\mathcal{P}_N(A))_3 < \infty.$
\end{thm}

\noindent Similarly, part $(ii)$ of Theorem \ref{main_thm_2} is
implied by the following
\begin{thm}[Breakdown of Universality for the Sherrington-Kirkpatrick
model]\label{(c):thm1} Let $\nu(r)=r^2$. For every bounded Borel set
A
\begin{gather}
\label{1st}
      \lim \mathbb{E} (\mathcal{P}_N(A))_1 =
      \mu(A), \quad \mathrm{P}\mbox{-a.s.}
\end{gather}
Moreover, if $\ \limsup  \frac{M(N)}{N} = \varepsilon < \frac{1}{8
\log 2}$ then $\mathrm{P}$-a.s.

(i) $\limsup  \mathbb{E} (\mathcal{P}_N(A))_2 =
\frac{(\mu(A))^2}{\sqrt{1 - 4 \varepsilon \log 2 }};$

(ii) $\limsup  \mathbb{E} (\mathcal{P}_N(A))_3 < \infty.$
\end{thm}
\begin{rem}
Condition $ \varepsilon < \frac{1}{8 \log 2}$ in not optimal and
could be improved. The reason for such a choice is that for
$\varepsilon < \frac{1}{8 \log 2}$ the third moment estimate is
quite simple.
\end{rem}

We will prove in detail Theorem \ref{(c):thm1} but omit the proof of
Theorem \ref{(c):thm2} as it is very similar and much simpler.

\begin{proof}[Proof of Theorem~\ref{(c):thm1}]
We successively prove the statement on the first, second, and third
moment.

\bigskip
\noindent 1.  {\bf  First moment estimate.}

Since all the random variables $H'_N(\sigma), \sigma \in S_N$, are
identically distributed then
\begin{equation} \mathbb{E}(\mathcal{P}_N(A))_1 = |X|
\mathbb{P}(H'_N(\sigma) \in A). \end{equation} We saw in the proof
of Theorem \ref{(ce):thm2} that for $M(N)$ satisfying $\log N =
o(2^M)$, $|X|=2^M(1+o(1))\,\,\mathrm{P}$-a.s. Combined with
\eqref{(a):common_term}, the definition of $a_N$, and the
dominated convergence theorem this fact implies that
$\mathrm{P}$-a.s.
\begin{align} \mathbb{E}(\mathcal{P}_N(A))_1 = 2^M(1+o(1))
\frac{1}{2^M} \frac{1}{\sqrt{\pi}} \int \limits_{A} e^{-x^2 b_N^2/2
- a_N b_N x} dx = \mu(A)(1+o(1)).
\end{align}
Hence \eqref{1st} is proven.

\comment{If $M(N) \leq \varepsilon N$ then to prove that the second
factorial moment $\mathbb{E}(\mathcal{P}_N(A))_2$ is uniformly in
$N$ bounded it suffices to prove this statement for the sequence
$M(N)=\varepsilon N.$ Indeed, assume there are two sequences such
that $M_1(N) \leq M_2(N)$ for all $N \in \mathbb{N}.$ We build on a
common probability space two families of i.i.d. Bernoulli random
variables $(X^1_{\sigma})_{\sigma \in S_N}$ and
$(X^2_{\sigma})_{\sigma \in S_N}$ such that $\mathbb{P}(X_{\sigma}^1
= 1) = \frac{2^{M_{1}}}{2^N}$ and $\mathbb{P}(X_{\sigma}^1 = 1) =
\frac{2^{M_{2}}}{2^N}.$ Since $M_1(N) \leq M_2(N)$ we can couple for
every $\sigma \in S_N$ the random variables $X_{\sigma}^1$ and
$X_{\sigma}^2$ in the following way: if $X_{\sigma}^1=0$ then
$X_{\sigma}^2=0,$ if $X_{\sigma}^1=1$ then $X_{\sigma}^1=1$ with
probability $2^{M_2-M_1}$ and $X_{\sigma}^1=0$ with probability
$1-2^{M_2-M_1}.$ From construction, the random cloud $X_1 =
\{\sigma\in S_N: X^1_{\sigma}=1\}$ contains the random cloud
$X_2=\{\sigma\in S_N: X^2_{\sigma}=1\}$ corresponding to the
sequence $M_2(N).$ As a consequence, for every $\ell \in \mathbb{N}$
the $\ell^{\textrm{th}}$ factorial moment corresponding to the
sequence $M_1(N)$ is greater than the one corresponding to the
sequence $M_2(N).$}


\bigskip
\noindent 2. {\bf  Second moment estimate.} Next assume that
$\limsup \frac{M(N)}{N} = \varepsilon$. We now want to calculate
$\limsup \mathbb{E}(\mathcal{P}_N(A))_2$. For this we rewrite the
second factorial moment as
\begin{equation}\label{(c)_moment}
\mathbb{E}(\mathcal{P}_N(A))_2 = \sum \limits_{
\sigma^{1},\sigma^{2}} \mathbb{P}\left( H'_{N}(\sigma^{1}) \in
A,H'_{N}(\sigma^{2}) \in A\right ),
\end{equation}
where the summation is over all pairs of distinct configurations
$(\sigma^{1},\sigma^{2}) \in X^{2}.$ We then split the set $S_N^{
2}$ into four non-intersecting subsets and calculate the
contributions from these subsets separately (these calculations are
similar to those of Theorem \ref{lem2}).

\medskip

(1) We begin by calculating the contribution from the set
$\mathcal{S}_1^X,$ where
\begin{align}
& \mathcal{S}_1 = \left\{ (\sigma^1,\sigma^2) \in S_N^2 : \
|R(\sigma^1,\sigma^2)| \leq \tau_N \right\}
\end{align}
and where the sequence $\tau_N$ is chosen in such a way that $N
\tau_N^4 \to 0$ and $e^{-N \tau_N^2}$ decays faster than any
polynomial: this can be achieved by choosing e.g. $\tau_N =
\frac{1}{N^{1/4} \log N}.$

\noindent First, we get from \eqref{(a):common_term} for $\ell = 2$
\begin{align}\label{(c):common_term}
\mathbb{P} \bigl( H'_{N}&(\sigma^{1}) \in A,H'_{N}(\sigma^{2}) \in
A\bigr) = \frac{1}{2^{2M}}\frac{e^{\frac{2b_{12}}{1+b_{12}} (M
\log 2 - 1/2 \log M)}}{\sqrt{1-b_{12}^2}}  \nonumber \\
& \times \frac{1}{(\sqrt{\pi})^2}\int \limits_{A}\int \limits_{A}
e^{-(\vec{x}, B^{-1} \vec{x}) b_N^2/2 - a_N b_N (\vec{x}, B^{-1}
\vec{1})} dx_1 dx_2,
\end{align}
where $b_{12} = \cov(H'_{N}(\sigma^{1}),H'_{N}(\sigma^{2})) =
R_{12}^2.$ If $(\sigma^1,\sigma^2) \in \mathcal{S}_1$ then $b_{12} =
o(1)$ and we get from \eqref{(a):int_part} that uniformly in
$(\sigma^1, \sigma^2) \in \mathcal{S}_1$ the second line in
\eqref{(c):common_term} is just $(\mu(A))^2(1+o(1)).$

\par Next, we represent the set $ \mathcal{S}_1$ as a union of
sets $V_2(R_{N,k})$ with $R_{N,k} = 1 - \frac{2k}{N} \in K_N \cap
[0, \tau_N].$ Applying Theorem~\ref{(ce):thm3} $(i)$ we get that
\begin{equation} |V_2(R_{N,k}) \cap X^{2}| = 2 \sqrt{\frac{2}{\pi N}}
2^{2M} e^{-N\mathcal{J}(R_{N,k})} (1+o(1)). \end{equation}

\noindent Therefore, up to a multiplicative term of the form
$1+o(1),$ the contribution from the set $\mathcal{S}_1^X$ to the sum
\eqref{(c)_moment} is
\begin{align}\label{(c):g:contr_S}
\sum \limits_{R_{N,k}} & \frac{1}{2^{2M}}
e^{\frac{2b_{12}}{1+b_{12}}(M\log 2 - 1/2 \log M)} (\mu(A))^2 2
\sqrt{\frac{2}{\pi N}} 2^{2M}
e^{-N\mathcal{J}(R_{N,k})} \nonumber \\
& = 2 \sqrt{\frac{2}{\pi N}}(\mu(A))^2 \sum \limits_{R_{N,k}} e^{-N
\mathcal{J}(R_{N,k})} e^{\frac{2b_{12}}{1+b_{12}}(M\log 2 - 1/2 \log
M)}.
\end{align}

Since $b_{12}(R_{N,k}) = R_{N,k}^2 > 0$ the last sum is monotone in
$M.$ As we will see below due to this fact it is sufficient to
calculate the upper limit of $\mathbb{E}(\mathcal{P}(A))_2$ for
sequences of the form $M(N) = \varepsilon N$ with $\varepsilon \in
(0,\frac{1}{8 \log 2}).$ Thus let $M = \varepsilon N, \varepsilon
\in (0,\frac{1}{8 \log 2}).$ We obtain
\begin{align}\label{(c):g:exp}
\frac{2b_{12}}{1+b_{12}} &(M \log 2 - 1/2 \log M) = 2 \varepsilon N
R_{N,k}^2 \log 2 +  O(N R_{N,k}^4).
\end{align}

\noindent The choice of $\tau_N$ guarantees that for $R_{N,k} \leq
\tau_N$ the last term in the rhs of \eqref{(c):g:exp} is of order
$o(1).$ Moreover, for $x<1$ we observe that
$\mathcal{J}(x)=\frac{1}{2}x^2 + O(x^4).$ Thus
\begin{equation}\label{(c):g:exp1}
N \mathcal{J}(R_{N,k}) = \frac{NR_{N,k}^2}{2} + o(1).
\end{equation}
Using \eqref{(c):g:exp} and \eqref{(c):g:exp1} the sum
\eqref{(c):g:contr_S} becomes
\begin{equation}\label{(c):g:contr_S(1)}
2 \sqrt{\frac{2}{\pi N}} (\mu(A))^2\sum \limits_{R_{N,k}} e^{-
\frac{1}{2}N R_{N,k}^2(1 - 4 \varepsilon \log 2)} (1+ o(1)),
\end{equation}
where the summation is over $R_{N,k}\in K_N \cap [0, \tau_N].$
Introducing new variables $y_{N,k} = \sqrt{N} R_{N,k}$ we further
rewrite \eqref{(c):g:contr_S(1)} as
\begin{equation}\label{(c):g:contr_S(2)}
\sqrt{\frac{2}{\pi}} \sum \limits_{y_{N,k}} \frac{2}{\sqrt{N}} e^{-
\frac{y_N^2}{2}(1 - 4 \varepsilon \log 2)} (1+ o(1)).
\end{equation}
It is not difficult to see that for $\varepsilon < \frac{1}{4 \log
2}$ the sum in \eqref{(c):g:contr_S(2)} converges to the integral
\begin{equation}
\int \limits_{0}^{\infty} e^{-\frac{y^2}{2}(1 - 4 \varepsilon \log
2)} dy = \frac{1}{\sqrt{1 - 4 \varepsilon \log2}}
\sqrt{\frac{\pi}{2}}.
\end{equation}
Therefore the contribution from the set $S_1^X$ is
$\frac{(\mu(A))^2}{\sqrt{1- 4 \varepsilon \log 2}}(1+o(1)).$

We can extend this result to the case when $M(N) = \varepsilon N +
o(N)$ with $\varepsilon < \frac{1}{8 \log 2}.$ Indeed, assume that
$\varepsilon_n \uparrow \varepsilon$ as $n \to \infty.$ Then using
the above calculation for $M = \varepsilon_n N$ together with the
monotonicity argument we have that for all $n \geq 1$
\begin{align}
\mathbb{E}(\mathcal{P}_N(A))_2 \geq \frac{(\mu(A))^2}{\sqrt{1- 4
\varepsilon_n \log 2}}.
\end{align}
Taking the limit in $n$ we obtain a lower bound. With exactly the
same argument we prove the corresponding upper bound.

\medskip

(2) Next, we estimate the contribution from the set
$\mathcal{S}_2^X,$ where
\begin{align}
\mathcal{S}_2 = \bigl\{ (\sigma^1,\sigma^2) \in S_N^2: \
|R(\sigma^1,\sigma^2)| > \tau_N \mbox{ and } R(\sigma^1,\sigma^2)
\mbox{ satisfies \eqref{(ce):thm3:1}}\bigr\}.
\end{align}

\par Since the set $A$ is bounded an elementary computation yields
that for a constant $C=C(A),$ uniformly in $(\sigma^1,\sigma^2) \in
\mathcal{S}_2,$
\begin{align}\label{(c):g:int2}
& \int \limits_{A} \int \limits_{A} e^{-(\vec{x}, B^{-1}\vec{x}) \
b_N^2/2 - a_N b_N (\vec{x}, B^{-1}\vec{1})} dx_1 dx_2 \leq C.
\end{align}

\par For fixed $N$ we let $R_{N}^{1}$  to be the largest value of
the overlap satisfying condition \eqref{(ce):thm3:1}. Then
representing $ \mathcal{S}_2 $ as a union of sets $V_2 (R_{N,k}),$
where $R_{N,k} = 1 - \frac{2k}{N} \in K_N \cap [\tau_N, R_N^1]$, and
using Theorem~\ref{(ce):thm3}~$(i)$, we conclude that the
contribution from the set $\mathcal{S}_2^X$ is, up to a
multiplicative term of the form $1+o(1),$ bounded~by
\begin{align}\label{(c):contr:S(2)}
\sum \limits_{R_{N,k}} 2 \sqrt{\frac{2}{\pi N(1 - R_{N,k}^2)}} \
\frac{C}{\sqrt{1-b_{12}^2}} \ e^{ -N \mathcal{J}(R_{N,k}) +
\frac{2b_{12}}{1+b_{12}} M \log 2}.
\end{align}

\par This quantity is monotone in $M.$ As a consequence,
to show that it is negligible in the limit $N \to \infty$ it
suffices to show this fact under the assumption $M(N)=\varepsilon
N.$ Thus, letting $M = \varepsilon N$ and using that $\mathcal{J}(x)
\geq \frac{1}{2} x^2$, we can bound the exponent in
\eqref{(c):contr:S(2)}~by
\begin{align}
-N \left( \mathcal{J}(R_{N,k}) - \frac{2b_{12}}{1+b_{12}}
\varepsilon \log 2 \right) \leq -\frac{1}{2} N x^2 ( 1 - 4
\varepsilon \log 2).
\end{align}

\noindent Since $1-b_{12}^2 \geq 1/N$ and since the number of terms
in \eqref{(c):contr:S(2)} is at most $N,$ we can further bound it,
for some positive constant $C>0$, by
\begin{align}
C N \ 2 \sqrt{\frac{2}{\pi}}e^{ - \frac{1}{2}N \tau_{N}^2(1 - 4
\varepsilon \log 2 )}.
\end{align}
Since $\varepsilon < \frac{1}{8 \log 2}$ the contribution from the
set $\mathcal{S}_2^X$ is negligible by definition of $\tau_N.$

\medskip

$(3)$ We now analyze the contribution from the set
$\mathcal{S}_3^X,$ where
\begin{align}
\mathcal{S}_3 = \bigl\{ (\sigma^1,\sigma^2) \in S_N^2:
R(\sigma^1,\sigma^2) \mbox{ satisfies \eqref{(ce):thm3:2}} \bigr\}.
\end{align}

\noindent Let $R_N^{2}$ be the largest overlap value  satisfying
condition \eqref{(ce):thm3:2}. To bound the contribution from the
set $\mathcal{S}_3$ we represent it as a union of sets $V_2
({R_{N,k}}),$ where $R_{N,k} = 1 - \frac{2k}{N}$ runs over the set
$K_N \cap [R_N^1, R_N^2].$ Proceeding as in~$(2)$ and using
Theorem~\ref{(ce):thm3} $(ii)$, we bound it by
\begin{align}\label{(c):contr:S(3)}
\sum \limits_{R_{N,k}} 2 \sqrt{\frac{2}{\pi N}} N^{c}
\frac{C}{\sqrt{1 - b_{12}^2}} e^{\frac{2b_{12}}{1+b_{12}} M \log 2 -
N \mathcal{J}(R_{N,k})}.
\end{align}
Again, the sum is monotone in $M$ and thus it is enough to bound it
for $M = \varepsilon N.$

\noindent For $R_{N,k}$ satisfying condition \eqref{(ce):thm3:2} we
can bound the exponent in \eqref{(c):contr:S(3)} as follows:
\begin{align}
\frac{2b_{12}}{1+b_{12}} & \varepsilon N \log 2 - N
\mathcal{J}(R_{N,k}) \leq \frac{2b_{12}}{1+b_{12}} \varepsilon N
\log 2 - N \varepsilon \log 2 - c_1\log N \nonumber \\ & = -
\frac{1-b_{12}}{1+b_{12}} \varepsilon N \log 2 - c_1 \log N \leq - C
N \varepsilon \log 2 - c_1\log N,
\end{align} where $C$ is
some positive constant.

\par Since $1 - b_{12}^2 \geq 1/N$ and since there are at most $N$ terms in the sum
\eqref{(c):contr:S(3)}, we can bound the latter by $ N^{c+1-c_1}
\exp\bigl( - C N \varepsilon \log 2\bigr).$ Therefore,
$\mathrm{P}$-almost surely the contribution from the set
$\mathcal{S}_3^X$ is negligible as $N \to \infty.$

\medskip

$(4)$ To finish the second moment estimate it remains to treat the
set
\begin{equation}
\mathcal{S}_4 = \bigl\{ (\sigma^1,\sigma^2) \in S_N^2:
R(\sigma^1,\sigma^2) \mbox{ satisfies } \eqref{(ce):thm3:3}\bigr\}.
\end{equation}
But by Theorem \ref{(ce):thm3} $(iii)$, the set $\mathcal{S}_4^X$ is
$\mathrm{P}$-almost surely empty. This finishes the proof of
assertion (i) of  Theorem~\ref{(c):thm1}.

\bigskip
\noindent 3.  {\bf  Third moment estimate.} To analyze the third
factorial moment we  use that, by formula \eqref{(a):gen_mom} and
definition \eqref{(ce):W}, it can be written as
\begin{align}\label{(b):p:2:m:3:moment}
\mathbb{E}\bigl( & \mathcal{P}_N(A)\bigr)_3 =
\sum\limits_{R_{12},R_{23},R_{13}} \frac{|W_{N,3}^X|}{2^{3M}
\sqrt{\det {B}}} \, e^{(3-(\vec{1}, B^{-1}\vec{1}))(M \log 2 - 1/2
\log M)} \nonumber \\ & \times \int \limits_{A} \int \limits_{A}
\int \limits_{A} e^{-{(\vec{x}, B^{-1}\vec{x})}b_N^2/{2} - a_N b_N
(\vec{x}, B^{-1}\vec{1}) } d \vec{x},
\end{align}
where $B = B(\sigma^1, \sigma^2, \sigma^3)$, the covariance matrix
of the vector $(H_N(\sigma^1), H_N(\sigma^2), H_N(\sigma^3))$, is
the matrix with elements $b_{ij} = R_{ij}^2, 1 \le i,j \le 3,$ and
where the summation runs over all triplets of overlaps
$(R_{12},R_{23},R_{31}) \in K_N^3.$ To estimate this sum we rely on
three auxiliary lemmas whose proofs we skip since they are simple.
\begin{lem}\label{(c):lem1}
If $\varepsilon < \frac{1}{8 \log 2}$ then
\begin{equation}
\limsup \limits_{N \to \infty} \max \limits_{\sigma^1,\sigma^2,
\sigma^3 \in X} \int \limits_{A} \int \limits_{A} \int \limits_{A}
e^{-{(\vec{x}, B^{-1}\vec{x})}b_N^2/{2} - a_N b_N (\vec{x},
B^{-1}\vec{1}) } d\vec{x} < \infty, \quad \mathrm{P}\mbox{-a.s.}
\end{equation}
\end{lem}

\begin{lem}\label{(c):lem2}
$\mathrm{P}$-almost surely for all configurations $\sigma^1,
\sigma^2,\sigma^3 \in X$
\begin{equation}\label{b:p:2:ini:ineq}
3 - (\vec{1}, B^{-1} \vec{1}) \leq 2(R_{12}^2 + R_{23}^2 +
R_{31}^2).
\end{equation}
\end{lem}

\begin{lem}\label{(c):lem3}
For all $\sigma^1, \sigma^2, \sigma^3 \in S_N$
\begin{equation}\label{(b):p:2:m:3:J}
\mathcal{J}^{(2)}(R_{12}, R_{23}, R_{31}) \ge \frac{1}{4} \Bigl(
R_{12}^2 + R_{23}^2 + R_{31}^2 \Bigr).
\end{equation}
\end{lem}

\par We are now ready to estimate sum \eqref{(b):p:2:m:3:moment}.
In the same spirit as for the second moment calculation, we will
split \eqref{(b):p:2:m:3:moment} into four parts and show that
every each of them is bounded. Moreover, by monotonicity argument
similar to that used in calculation of the second moment we can
restrict our attention to the case when $M = \varepsilon N.$

\medskip

$(1)$ We first calculate the contribution to
\eqref{(b):p:2:m:3:moment} coming from the set
\begin{align}
\mathcal{S}_1 = \bigl\{ (\sigma^1,\sigma^2, \sigma^3) \in S_N^3: \
\max \bigl\{ |R_{12}|,|R_{23}|,|R_{31}| \bigr\} \le \tau_N\bigr\},
\end{align}
where $\tau_N = \frac{1}{N^{1/4} \log N}.$ Using Theorem
\ref{(ce):thm6} $(i)$ and Lemma \ref{(c):lem2}, we obtain that the
contribution from $\mathcal{S}_1$ is at most of  order
\begin{align}\label{(b):p:2:m:3:moment:1}
\frac{1}{N^{3/2}} \sum \limits_{R_{12},R_{23},R_{31}}
e^{-N\mathcal{J}^{(2)}(R_{12},R_{23},R_{31}) +
2(R_{12}^2+R_{23}^2+R_{31}^2)M\log 2},
\end{align}
where the summation is over $R_{12},R_{23},R_{31} \in K_N \cap
[-\tau_N, \tau_N].$ Expanding in Taylor series we obtain that for
$|R_{12}|,|R_{23}|,|R_{31}| \leq \tau_N$
\begin{equation}
\mathcal{J}^{(2)}(R_{12}, R_{23}, R_{31}) = \frac{1}{2}\bigl(
R_{12}^2 + R_{23}^2  + R_{31}^2 \bigr) + O(\tau_N^4),
\end{equation}
and thus  \eqref{(b):p:2:m:3:moment:1} is bounded by
\begin{align}
\frac{1}{N^{3/2}} \sum \limits_{R_{12},R_{23},R_{31}} &
e^{-\frac{1}{2} N( R_{12}^2 + R_{23}^2  + R_{31}^2) +
2\varepsilon N (R_{12}^2+R_{23}^2+R_{31}^2)\log 2} \nonumber \\
& = \left( \frac{1}{\sqrt{N}} \sum \limits_{R_N \in K_N \cap
[-\tau_N,\tau_N]} e^{-\frac{1}{2}NR^2(1 - 4 \varepsilon \log 2)}
\right)^3 < \infty.
\end{align}

\medskip

(2) We next calculate the contribution from the set
\begin{align}
\mathcal{S}_2 = \bigl\{ (\sigma^1, \sigma^2, \sigma^3) \in S_N^3: \,
& R_{12},R_{23},R_{21} \mbox{ satisfy \eqref{(ce):thm6:1}} \nonumber
\\ & \mbox{and } \max \bigl\{|R_{12}|,|R_{23}|,|R_{31}| \bigr \} >
\tau_N \bigr\}.
\end{align}
Without loss of generality we can assume that $|R_{12}| > \tau_N.$
Then, using Theorem \ref{(ce):thm6} $(i)$ and Lemma \ref{(c):lem3},
the contribution from this set is at most of order
\begin{align}\label{(b):p:2:m:3:moment:2}
\sum \limits_{R_{12},R_{23},R_{31}}
\frac{e^{-N\mathcal{J}^{(2)}(R_{12},R_{23},R_{31})}}{N^{3/2}
P(R_{12},R_{23},R_{31})} e^{2(R_{12}^2 + R_{23}^2 + R_{31}^2) M \log
2},
\end{align}
where $P(x,y,z)=(1+ x + y + z)(1+x - y - z)(1 - x + y - z)(1 - x - y
+ z)$ and where the sum is over the triplets $(R_{12},R_{23},R_{31})
\in K_N^3$ satisfying \eqref{(ce):thm6:1} and $|R_{12}| > \tau_N.$
Then, from Lemma \ref{(c):lem2}, we further get that the sum
\eqref{(b):p:2:m:3:moment:2} is bounded by
\begin{align}
\sum \limits_{|R_{12}| > \tau_N} & \frac{e^{-\frac{1}{4} N(R_{12}^2
+ R_{23}^2 + R_{31}^2)}}{N^{3/2} P(R_{12},R_{23},R_{31})}
e^{2(R_{12}^2 + R_{23}^2 + R_{31}^2) M \log 2} \nonumber \\ & = \sum
\limits_{|R_{12}| > \tau_N} \frac{e^{-\frac{1}{4} N (1 - 8
\varepsilon \log 2)( R_{12}^2 + R_{23}^2 + R_{31}^2)}}{N^{3/2}
P(R_{12},R_{23},R_{31})}
\end{align}
Since the number of terms in the sum is at most $(N+1)^3$, and since
$e^{-\frac{1}{4} (1 - 8 \varepsilon \log 2) N \tau_N^2}$ decreases
faster than any polynomial it follows that the sum is of order
$o(1).$

\medskip

(3) We now turn to the contribution from the set
\begin{align}
\mathcal{S}_3 = \bigl\{ (\sigma^1, & \sigma^2, \sigma^3) \in
S_N^3: \, R_{12},R_{23},R_{31} \mbox{ satisfy
\eqref{(ce):thm6:2}}\bigr\}.
\end{align}
By Lemma \ref{(ce):thm6} $(ii)$ the contribution from this set is at
most of order
\begin{align}\label{(b):p:2:m:3:moment:3}
\sum \limits_{R_{12},R_{23},R_{31}} & \frac{N^c \, \mathrm{E}
|W_{N,3}^X|}{2^{3M} \sqrt{\det B}} e^{(3-(\vec{1},B^{-1}\vec{1}))(M
\log 2 - 1/2 \log M)},
\end{align}
where the summation is over the triplets $(R_{12},R_{23},R_{31}) \in
K_N$ satisfying \eqref{(ce):thm6:2}. Since $\mathrm{E} |W_{N,3}^X|$ is
of order
\begin{equation}
\frac{2^{3M} e^{-N \mathcal{J}^{(2)}(R_{12},R_{23},R_{31})}}{N^{3/2}
P(R_{12},R_{23},R_{31})}
\end{equation}
we obtain, using Lemmas \ref{(c):lem2} and \ref{(c):lem3}, that the
sum \eqref{(b):p:2:m:3:moment:3} is bounded by
\begin{equation}
\sum \limits_{R_{12},R_{23},R_{31}} \frac{N^c e^{-\frac{1}{4}
N(R_{12}^2 + R_{23}^2 + R_{31}^2)}}{N^{3/2} P(R_{12},R_{23},R_{31})}
e^{2\varepsilon N \log 2 (R_{12}^2 + R_{23}^2 + R_{31}^2)}.
\end{equation}
It is easy to show that the triplet $(R_{12},R_{23},R_{31})$ that
satisfy \eqref{(ce):thm6:2} must satisfy either $|R_{23}| > \tau_N$
or $|R_{31}| > \tau_N.$ Therefore we can further bound the
contribution from the set $\mathcal{S}_3$ by the sum
\begin{equation}
\sum \limits_{|R_{23}| > \tau_N \mbox{ or } |R_{31}| > \tau_N}
\frac{N^c e^{-\frac{1}{4} N(R_{12}^2 + R_{23}^2 +
R_{31}^2)}}{N^{3/2} P(R_{12},R_{23},R_{31})} e^{2\varepsilon N \log
2 (R_{12}^2 + R_{23}^2 + R_{31}^2)},
\end{equation}
which is $o(1)$ by the same argument as in part $(2).$

\medskip

(4) To finish the estimate of the third factorial moment we have to
estimate the contribution to \eqref{(b):p:2:m:3:moment} coming from
the set
\begin{equation}
\mathcal{S}_4 = \bigl\{ (\sigma^1, \sigma^2, \sigma^3) \in S_N^3: \,
R_{12},R_{23},R_{21} \mbox{ satisfy \eqref{(ce):thm6:3}}\bigr\}.
\end{equation}
By Theorem \ref{(ce):thm6} $(iii)$ the set $\mathcal{S}_4^X$ is
$\mathrm{P}$-a.s. empty and therefore its contribution is
$\mathrm{P}$-a.s. zero. This finishes the proof of assertion $(ii)$
of Theorem \ref{(c):thm1}. The proof of  Theorem~\ref{(c):thm1} is
now complete.
\end{proof}

\bigskip

For comparison with the cases $\nu(r)= r^p$ for $p=1$ and $p=2$, we
give here the asymptotic behavior of the first three factorial
moments for the case $\nu(0)= \nu'(0)=0$, i.e for instance for the
case $\nu(r)= r^p$ of pure $p$-spins when $p \ge 3.$ This behavior
is compatible with a Poisson convergence theorem.

\begin{thm}\label{(c):thm3}
Assume that $\nu(0)= \nu'(0)=0$.  For every bounded Borel set A
\begin{equation}
\lim \mathbb{E} (\mathcal{P}_N(A))_1 = \mu(A), \quad
\mathrm{P}\mbox{-a.s.}
\end{equation}
Moreover, if $\ \limsup \frac{M(N)}{N} < \frac{1}{8 \log 2}$ then
$\mathrm{P}$-a.s.

(i)  $\lim \mathbb{E} (\mathcal{P}_N(A))_2 = (\mu(A))^2;$

(ii) $\lim \mathbb{E} (\mathcal{P}_N(A))_3 = (\mu(A))^3 < \infty.$
\end{thm}
We do not include a proof of this last statement, which again
follows the same strategy as the proof of Theorem \ref{(c):thm1}.

\bigskip

\section{Universality for Non-Gaussian Hamiltonians}\label{NG}

\par In this section we extend the results of the previous sections
to the case of  non-Gaussian Hamiltonians. We are able to make
this extension only for the pure $p$-spin models, i.e. $\nu(r)=
r^p$. In this case we recall that the Hamiltonian is defined~as
\begin{equation}\label{(ng):Hamiltonian}
H_N(\sigma)=\frac{1}{\sqrt{N}} H_{N,p} = \frac{1}{\sqrt{N^p}} \sum
\limits_{1 \leq i_1, \dots, i_p \leq N} g_{i_1,\dots,i_p}
\sigma_{i_1} \dots \sigma_{i_p}.
\end{equation}
Our assumptions on the random variables $(g_{i_1,\dots,i_p})_{1
\leq i_1, \dots, i_p \leq N}$ in \eqref{(ng):Hamiltonian} are the
same as were made in \cite{BCMN1} and \cite{BCMN2} for the number
partitioning problem. That is, we assume that their distribution
function admits a density $\rho(x)$ that satisfies the following
conditions:
\begin{enumerate}
  \item $\rho(x)$ is even;
  \item $\int x^2 \rho(x) dx = 1;$
  \item for some $\epsilon > 0$
  \begin{equation}\label{(intro):rho}
  \int \limits_{-\infty}^{\infty} \rho(x)^{1+\epsilon} dx <
  \infty;
  \end{equation}
  \item $\rho(x)$ has a Fourier transform that is analytic in some
neighborhood of zero. We write
\begin{equation}\label{(ng):F}
- \log \hat{\rho} (z) = \frac{1}{2}(2\pi)^2 z^2 + c_4(2\pi)^4 z^4 +
O(|z|^6).
\end{equation}
Note that the inequality $\mathbb{E}(X^4) \ge \mathbb{E}(X^2)^2$
implies that necessarily $c_4 < \frac{1}{12}.$
\end{enumerate}

\par Under these assumptions we will show, using the method introduced
by C. Borgs, J. Chayes, S. Mertens and C. Nair in \cite{BCMN2}, that
Theorems \ref{main_thm} and \ref{main_thm_2} still hold.
%

\subsection{Proof of Universality.}\label{(ng.a)}

In this subsection we fix $p \ge 1$ and prove the analog of
Theorem \ref{main_thm} in the non-Gaussian case assuming that the
Hamiltonian is given by \eqref{(ng):Hamiltonian}, and that the
random variables $(g_{i_1,\dots,i_p})_{1 \leq i_1, \dots, i_p \leq
N}$ satisfy conditions $(1)-(4)$ above.

\begin{thm}[Universality in the Non-Gaussian case]
\label{(ng):main_thm} Assume $M(N) = o(\sqrt{N})$ for $p=1$ and
$M=o(N)$ for $p \geq 2.$ Then $\mathrm{P}$-almost surely the
sequence of point processes $\mathcal{P}_N$ converges weakly to a
Poisson point process $\mathcal{P}$ on $\mathbb{R}$ with intensity
given by
\begin{equation}
\mu(dt) = \frac{1}{\sqrt{\pi}} e^{-t \sqrt{2 \log 2}} dt.
\end{equation}
\end{thm}

To prove Theorem \ref{(ng):main_thm} we essentially prove a local
limit theorem. More precisely, for any fixed $N$ let us introduce
the Gaussian process $Z_N$ on $S_N$ that has the same mean and
covariance matrix as the process $H'_{N}(\sigma)$ defined in
\eqref{(intro):Hamiltonian:normalized}. We will prove in Theorem
\ref{(ng):thm1} that $\mathrm{P}$-a.s., for all sequences
$(\sigma^1,\dots,\sigma^\ell) \in X^{\ell}$, the joint density of
the random variables $H'_N(\sigma^1),\dots,H'_N(\sigma^\ell)$ is
well approximated by the joint density of
$Z_N(\sigma^1),\dots,Z_N(\sigma^\ell).$

\begin{thm}\label{(ng):thm1}
Assume $M(N)= o(\sqrt{N})$ for $p=1$ and $M=o(N)$ for $p \geq 2.$
Then $\mathrm{P}$-almost surely for every $\ell \ge 1$ and every
bounded Borel set $A$ there exists $c>0$ such that uniformly in
$(\sigma^1,\dots,\sigma^\ell) \in X^{\ell}$
\begin{align}\label{(ng):approx:denst}
\mathbb{P}( H'_N(\sigma^j) & \in A,j=1,\dots,\ell) = \mathbb{P}(
Z_N(\sigma^j) \in A,j=1,\dots,\ell) \nonumber \\ & \times
\Bigl(1+O(R_{\max}(\sigma^1,\dots,\sigma^\ell)) + O
\Bigl(\frac{M^2}{N^p} \Bigr) \Bigr)  + O(e^{-cN^p}).
\end{align}
\end{thm}

Applying Theorem~\ref{(ng):thm1} together with
Theorem~\ref{main_thm} we get from formula \eqref{(a):gen_mom} that
\begin{align}
\mathbb{E} ( \mathcal{P}_N(A) )_{\ell} = (\mu(A))^\ell (1+o(1)) +
O(2^{M \ell} e^{-c N^p }) \to (\mu(A))^\ell,
\end{align}
which, by Lemma \ref{(a):gen_lem}, implies weak convergence of the
sequence of  point processes $\mathcal{P}_N$ to a Poisson point
process with intensity measure $\mu,$ thus implying Theorem
\ref{(ng):main_thm}. We therefore focus on the proof of Theorem
\ref{(ng):thm1}.

\begin{proof}

First, we obtain from the definition of $H'_N(\sigma)$ that
\begin{align}
\Bigl \{ H_N'(\sigma) \in (x,x+{\vartriangle}x) \Bigr\} = & \Bigl\{
\sum_{1 \leq i_1, \dots, i_p\leq N} g_{i_1,\dots,i_p}
\sigma_{i_1} \dots \sigma_{i_p} \in \nonumber \\
& ( a_N + x b_N, a_N + (x+ {\vartriangle}x)b_N) \sqrt{N^p} \Bigr\}.
\end{align}

\noindent Following \cite{BCMN2} we get an integral representation
of the indicator function
\begin{align}\label{(ng):int:repr}
\textbf{1}&_{H'_N(\sigma)\in (x,x+{\vartriangle}x)} \nonumber \\ & =
{\vartriangle}x \, b_N \sqrt{n} \int \limits_{-\infty}^{\infty}
\sinc \bigl(f{\vartriangle}x b_N\sqrt{n}\bigr) \, e^{ 2 \pi i f \sum
g_{i_1,\dots,i_p} \sigma_{i_1} \dots \sigma_{i_p} - \, 2 \pi i f
\alpha_N \sqrt{n}} df,
\end{align}
where, for brevity, we wrote $n=N^p, \alpha_N = a_N + b_N
(x+{\vartriangle}x/2), \sinc(x) = \frac{\sin(\pi x)}{\pi x}$, and
where the sum in the exponent runs over all possible sequences $1
\leq i_1,\dots,i_p\leq N.$

\par Changing the integration variable in \eqref{(ng):int:repr} from
$f$ to $-f$ and applying the resulting formula to the product of
indicator functions we arrive at the following representation
\begin{align}\label{(ng):int:repr1}
\prod \limits_{j=1}^{\ell} & \textbf{1}_{H'_N(\sigma^j)\in
(x_j,x_j+\vartriangle x_j)} = \prod \limits_{j=1}^{\ell}
{\vartriangle}x_j \, b_N \sqrt{n} \nonumber \\ & \times
\iiint_{-\infty}^{\infty} \prod\limits_{j=1}^{\ell} \sinc
\bigl(f_j{\vartriangle}x_j b_N \sqrt{n}\bigr) \, e^{ - 2 \pi i f_j
\sum g_{i_1 \dots i_p} \sigma_{i_1}^j \dots \sigma_{i_p}^j + \, 2
\pi i f_j \alpha_N^{(j)} \sqrt{n} }df_j,
\end{align}

\noindent where $\alpha_N^{(j)} = a_N + b_N
(x_j+{\vartriangle}x_j/2).$ Introducing the variables
\begin{equation}
v_{i_1,\dots,i_p} = \sum \limits_{j=1}^{\ell} f_j
\sigma^{j}_{i_1}\dots\sigma^{j}_{i_p},
\end{equation}
we rewrite the integral in the above formula as
\begin{equation}
\iiint_{-\infty}^{\infty} \prod \limits_{1 \leq i_1, \dots, i_p \leq
N} e^{-2 \pi i g_{i_1 \dots i_p} v_{i_1,\dots,i_p}}
\prod\limits_{j=1}^{\ell} \sinc \bigl(f_j{\vartriangle}x_j
\sqrt{n}\bigr) e^{2 \pi i f_j \alpha_N^{(j)} \sqrt{n}}df_j.
\end{equation}

\par To get an integral representation of the joint density
\begin{equation}\label{(ng):joint:den}
\mathbb{P}\bigl( H'_N(\sigma^1) \in
(x_1,x_1+dx_1),\dots,H'_N(\sigma^\ell) \in (x_\ell,x_\ell+dx_\ell)
\bigr)
\end{equation}
we have to take the expectation of \eqref{(ng):int:repr1} and then
let ${\vartriangle}x_j \to 0$ for all $j=1,2,\dots \ell$. As was
proved in \cite{BCMN1} (see Lemma 3.4), the exchange of expectation
and integration for $p=1$ is justified when the rank of the matrix
formed by the row vectors $\sigma^1,\dots,\sigma^\ell,$ is $\ell.$
To justify the exchange in our case we introduce an $\ell$ by
${N}^{p}$ matrix, $C_p(\sigma^1,\dots,\sigma^\ell),$ defined as
follows: for any given set of configurations
$\sigma^1,\dots,\sigma^\ell,$ the $j$-th row is composed of all
${N}^{p}$ products, $\sigma^j_{i_1}\sigma^j_{i_2}\dots
\sigma^j_{i_p}$ over all subsets $1\leq i_1, \dots, i_p \leq N.$ By
generalizing the arguments from \cite{BCMN1} the exchange can then
be justified provided that the rank of the matrix
$C_p(\sigma^1,\dots,\sigma^\ell)$ is $\ell.$ As we will see in Lemma
\ref{(ng):Cp:lem} below this holds true $\mathrm{P}$-almost surely
when $M=o(N).$

\par Given a vector
$\boldsymbol{\delta} \in \{-1,1\}^{\ell}$ let
$n_{\boldsymbol{\delta}}$ be the number of times the column vector
$\boldsymbol{\delta}$ appears in the matrix~$C_p:$
\begin{equation}
n_{{\boldsymbol{\delta}}}
=n_{\boldsymbol{\delta}}(\sigma^1,\dots,\sigma^\ell)= \left|
\left\{j\leq N^p:
(\sigma^1_{j},\dots,\sigma^{\ell}_{j})=\boldsymbol{\delta}
\right\}\right|.
\end{equation}
With this notation we have:

\begin{lem}\label{(ng):Cp:lem}
Suppose $M=o(N).$ Then there exists a sequence $\lambda_N = o(1)$
such that $\mathrm{P}$-almost surely for all collections
$(\sigma^1,\dots,\sigma^\ell) \in X^{\ell}$
\begin{equation}\label{(a):ng:Cp_n}
\displaystyle \max \limits_{\boldsymbol \delta \in \{-1,1\}^\ell
}\Bigl| n_{\boldsymbol \delta} - \frac{n}{2^\ell}\Bigr| \leq n
\lambda_N.
\end{equation}
\end{lem}

\begin{proof}

We first prove by induction that the following simple fact holds
true: if for a given sequence of configurations
$(\sigma^1,\dots,\sigma^\ell) \in S_N^{\ell}$ the matrix
$C_{1}(\sigma^1,\dots,\sigma^\ell)$ satisfies condition
\eqref{(a):ng:Cp_n} then necessarily the matrix
$C_{p}(\sigma^1,\dots,\sigma^\ell)$ satisfies \eqref{(a):ng:Cp_n}
for all $ p\geq 1.$

For $p=1$ there is nothing to prove. We now assume that the
statement is true for the matrix
$C_{p-1}(\sigma^1,\dots,\sigma^\ell)$
and prove it for $C_{p}(\sigma^1,\dots,\sigma^\ell).$

\par Let $\sigma_{\mu}, 1 \leq \mu \leq N$ denote the columns of the
matrix $C_1(\sigma^1,\dots,\sigma^\ell)$. For every column vector
$\sigma_{\mu}$ let us construct a matrix
$C_{p-1}^{\mu}=C_{p-1}^{\mu}(\sigma^1,\dots,\sigma^\ell)$ with
entries
\begin{equation}
(C_{p-1}^{\mu})_{ij} = (\sigma_{\mu})_j \; (C_{p-1})_{ij}.
\end{equation}
For future convenience let $n_{\boldsymbol \delta}^{\mu}$ denote the
variable $n_{\boldsymbol \delta}$ for the matrix $C_{p-1}^{\mu}.$
>From the inductive assumption it follows that for all $1 \leq \mu
\leq N$

\begin{equation}
\max \limits_{\boldsymbol \delta} \Bigl| n_{\boldsymbol
\delta}^{\mu} - \frac{N^{p-1}}{2^\ell}\Bigr| \leq N^{p-1} \lambda_N.
\end{equation}

\noindent Now note that the $\ell \times N^p$ matrix
$C_{p}(\sigma^1,\dots,\sigma^\ell)$ can be obtained by
concatenating $N$ matrices
$C_{p-1}^{\mu}(\sigma^1,\dots,\sigma^\ell)$ each of size $\ell
\times N^{p-1}.$ Therefore for any sequence of configurations
$(\sigma^1,\dots,\sigma^\ell)$ with matrix $C_{1}(
\sigma^1,\dots,\sigma^\ell )$ satisfying \eqref{(a):ng:Cp_n} we
have

\begin{align}
\Bigl| n_{\boldsymbol \delta} - \frac{N^p}{2^{\ell}} \Bigr| & \leq
\Bigl| n_{\boldsymbol \delta}^{\mu_1} - \frac{N^{p-1}}{2^{\ell}}
\Bigr| + \dots + \Bigl| n_{\boldsymbol \delta}^{\mu_N} -
\frac{N^{p-1}}{2^{\ell}}
\Bigr| \leq \nonumber \\
& \leq N^{p-1}\lambda_N + \dots + N^{p-1}\lambda_N = N^{p}\lambda_N
\end{align}
and the induction is complete.

\par To prove Lemma \ref{(ng):Cp:lem} it is thus enough
to demonstrate that $\mathrm{P}$-almost surely there are no
sequences $(\sigma^1,\dots,\sigma^\ell) \in X^{\ell}$ such that
$C_1(\sigma^1,\dots,\sigma^\ell)$ violates condition
\eqref{(a):ng:Cp_n}. Let us prove this by induction in $\ell.$

\par For $\ell = 1$ let us introduce the sets
\begin{equation}
\mathcal{T}_N = \left\{\sigma \in S: \max \limits_{\boldsymbol
\delta \in \{-1,1\}} \Bigl|n_{\boldsymbol \delta} -
\frac{N}{2}\Bigr| \leq N \lambda_N \right\}.
\end{equation}
Then by Chernoff bound
\begin{align}
|\mathcal{T}_N^c| = 2 \sum \limits_{ i \geq N \lambda_N }
\binom{N}{\frac{N}{2}+i} \leq 2^{N + 1} e^{- \frac{1}{2} N ((1 +
\lambda_N) \log(1 + \lambda_N) - \lambda_N) }.
\end{align}
Let us choose $\lambda_N = o(1)$ in such a way that $M=o(N
\lambda_N^2)$ and $\log N = o(N \lambda_N^2).$ Then using
\eqref{(ce):(aux):ldev:4} we obtain that
\begin{equation}
\sum \limits_{ \sigma \in \mathcal{T}_N^c} {\bf 1}_{\sigma} = 0,
\quad \mathrm{P}\mbox{-a.s.}
\end{equation}
which proves the statement for $\ell = 1.$ Now assume that
$\mathrm{P}$-a.s. for all sequences
$(\sigma^1,\dots,\sigma^{\ell-1}) \in X^{(\ell-1)}$ the matrix
$C_1(\sigma^1,\dots,\sigma^{\ell-1})$ satisfies condition
\eqref{(a):ng:Cp_n}. Since there is only a countable number of
sequences $(\sigma^1,\dots,\sigma^{\ell-1})$ we fix
$(\sigma^1,\dots,\sigma^{\ell-1})\in X^{(\ell -1)}$ and prove that
$\mathrm{P}$-almost surely there are no configurations $\sigma^\ell$
such that $C_1(\sigma^1,\dots,\sigma^\ell)$ violates
\eqref{(a):ng:Cp_n}.

\par Let $\boldsymbol \delta \in \{-1,1\}^\ell$ be given. Define
$\boldsymbol \delta_1(\boldsymbol \delta) \in \{-1,1\}^{\ell-1}$ as
$(\boldsymbol \delta_1)_i = (\boldsymbol \delta)_i$ for $1 \leq i
\leq \ell-1$ and also define $\boldsymbol \delta_2(\boldsymbol
\delta) \in \{-1,1\}$ as $\boldsymbol \delta_2 = (\boldsymbol
\delta)_\ell.$ Let us also introduce, for given $\boldsymbol
\delta_1 \in \{-1,1\}^{\ell-1}$, the set
\begin{equation}
N_{\boldsymbol \delta_1} = \{j \leq N:
(\sigma^1_j,\dots,\sigma^{\ell-1}_j)=\boldsymbol \delta_1\}.
\end{equation}
By the inductive assumption we conclude that for all $\boldsymbol
\delta_1 \in \{-1,1\}^{\ell-1}$
\begin{equation}\label{(ng):lem:Ndelta1}
\Bigl||N_{\boldsymbol \delta_1}| - \frac{N}{2^{\ell-1}}\Bigr| \leq N
\lambda_N.
\end{equation}
From the definition of $N_{\boldsymbol \delta}$ it is not hard to
see that for all $\boldsymbol \delta \in \{-1,1\}^{\ell}$
\begin{equation}
n_{\boldsymbol \delta} = \bigl|\{j  \in N_{\boldsymbol
\delta_1(\boldsymbol \delta)} : \sigma^\ell_j = \boldsymbol
\delta_2(\boldsymbol \delta) \}\bigr|.
\end{equation}
Using the above relation together with \eqref{(ce):(aux):ldev:3} and
the assumptions on $\lambda_N$ we get that $\mathrm{P}$-almost
surely
\begin{equation}\label{(ng):lem:Ndelta}
\Bigl| n_{\boldsymbol \delta} - \frac{|N_{\boldsymbol
\delta_1(\boldsymbol \delta)}|}{2}\Bigr| \leq \frac{N \lambda_N}{2}.
\end{equation}

\noindent Therefore from \eqref{(ng):lem:Ndelta1} and
\eqref{(ng):lem:Ndelta}
\begin{align}
\Bigl| n_{\boldsymbol \delta} - \frac{N}{2^{\ell}} \Bigr| \leq
\Bigl| n_{\boldsymbol \delta} - \frac{|N_{\boldsymbol
\delta_1(\boldsymbol \delta)}|}{2} \Bigr|  + \Bigl|
\frac{|N_{\boldsymbol \delta_1(\boldsymbol \delta)}|}{2} -
\frac{N}{2^\ell}\Bigr| \leq \frac{N \lambda_N}{2} + \frac{N
\lambda_N}{2} = N \lambda_N.
\end{align}
The induction is now complete and the lemma is proved.
\end{proof}

\medskip

\par Lemma \ref{(ng):Cp:lem} implies that $\mathrm{P}$-almost surely, for all
$(\sigma^1,\dots,\sigma^\ell) \in X^{\ell}$,
\begin{equation}\label{(a):n_min}
n_{\min}(\sigma^1,\dots,\sigma^\ell) = \min
\limits_{\boldsymbol{\delta} \in \{-1,1\}^{\ell}}
n_{\boldsymbol{\delta}} = \frac{n}{2^{\ell}}(1 + O(\lambda_N ) ),
\end{equation}

\noindent and hence, for sufficiently large $N$, the rank of the
matrix $C_p(\sigma^1,\dots, \sigma^{\ell})$ is $\ell.$ The exchange
of integration and expectation is thus justified.

\par Using once again Lemma~\ref{(ng):Cp:lem}, condition
\eqref{(intro):rho}, and the dominated convergence theorem we obtain
that the joint density is
\begin{align}\label{(ng):prob}
\mathbb{P}& \bigl( H'_N(\sigma^j) \in (x_j, x_j + dx_j) \mbox{ for }
j=1,\dots,\ell \bigr) = \prod \limits_{j=1}^{\ell} b_N \sqrt{n} \,
dx_j \nonumber \\ & \times \iiint_{-\infty}^{\infty} \prod \limits_{
1\leq i_1, \dots, i_p \leq N} \hat{\rho}(v_{i_1,\dots,i_p})
\prod\limits_{j=1}^{\ell} e^{2 \pi i f_j \alpha_N^{(j)} \sqrt{n}}
df_j,
\end{align}
where we redefined $\alpha_N^{(j)} = a_N + b_N x_j.$ We remark for
future use that $\alpha_N^{(j)} = O(a_N)$ for all $ 1 \leq j \leq
\ell.$

\par It is straightforward at this point to adapt the saddle point
analysis used in \cite{BCMN2} to calculate the integrals of such
type. The only difference is that instead of the matrix
$C_1(\sigma^1,\dots,\sigma^\ell)$ with rows formed by row vectors
$\sigma^1,\dots,\sigma^\ell$, we use  the matrix
$C_p(\sigma^1,\dots,\sigma^\ell).$ By analogy with Lemma 5.3 from
\cite{BCMN2} we first approximate the integral in \eqref{(ng):prob}
by an integral over a bounded domain, i.e. for some $c_1>0$
depending on $\mu_1>0$

\begin{align}\label{(a):ng:prob1}
\iiint_{-\infty}^{\infty} & \prod \limits_{ i_1, \dots, i_p}
\hat{\rho}(v_{i_1,\dots,i_p}) \prod\limits_{j=1}^{\ell} e^{2 \pi i
f_j \alpha_N^{(j)} \sqrt{n}} df_j = \nonumber \\ &
\iiint_{-\mu_1}^{\mu_1} \prod \limits_{ i_1, \dots, i_p}
\hat{\rho}(v_{i_1,\dots,i_p}) \prod\limits_{j=1}^{\ell} e^{2 \pi i
f_j \alpha_N^{(j)} \sqrt{n}} df_j + O(e^{-c_1 n_{\min}}) .
\end{align}

\noindent We next rewrite the integral in the  r.h.s.~of
\eqref{(a):ng:prob1} as
\begin{align}\label{(a):ng:prob2}
\iiint_{-\mu_1}^{\mu_1} e^{2 \pi i n \mathbf{f} \cdot
\boldsymbol{\alpha}} \prod \limits_{ \boldsymbol \delta \in
\{-1,1\}^\ell } \hat{\rho}({\mathbf f}\cdot \boldsymbol{\delta}
)^{n_{\boldsymbol{\delta}}} \prod\limits_{j=1}^{\ell} df_j,
\end{align}

\noindent where
$\boldsymbol{\alpha}=\bigl(\frac{\alpha_N^{(1)}}{\sqrt{n}},\dots,\frac{\alpha_N^{(\ell)}}{\sqrt{n}}\bigr)$,
$\boldsymbol f = (f_1,\dots,f_\ell)$, and where $\boldsymbol{\alpha}
\cdot \boldsymbol f = \alpha_1 f_1 + \dots + \alpha_\ell f_\ell$ is
the standard scalar product.

\par Using Lemma \ref{(ng):Cp:lem} again we can apply Lemma 5.4 from
\cite{BCMN2} to conclude that given $\mu_1$ there are constants
$c_1(\mu_1)>0$ and $\mu_2>0$ such that the following equality holds
whenever $\eta_1,\dots,\eta_\ell$ is a sequence of real numbers with
$\sum_j |\eta_j| \leq \mu_2$ and $\eta_j \alpha_N^{(j)} \geq 0$ for
all $j=1,\dots,\ell$

\begin{align}\label{(a):ng:prob3}
\iiint_{-\mu_1}^{\mu_1} & \prod \limits_{\boldsymbol \delta}
\hat{\rho} (\mathbf{f} \cdot \boldsymbol \delta)^{n_{\boldsymbol
\delta}} \prod \limits_{j=1}^{\ell} e^{2 \pi i n f_j
\alpha_N^{(j)}}df_j \\
& = \iiint_{-\mu_1}^{\mu_1} e^{2 \pi n ( i \mathbf{f} \cdot
\boldsymbol \alpha - \boldsymbol{\eta} \cdot \boldsymbol{\alpha}) }
\prod \limits_{\boldsymbol \delta} \hat{\rho} (\mathbf{f} \cdot
\boldsymbol \delta + i \boldsymbol \eta \cdot \boldsymbol \delta
)^{n_{\boldsymbol \delta}} \prod \limits_{j=1}^{\ell} df_j +
O(e^{-\frac{1}{2}c_1 n_{\min}}). \nonumber
\end{align}

\noindent The values of the shifts $\eta_1,\dots,\eta_\ell$ are
determined by the following system:
\begin{equation}\label{(a):ng:system:G}
\sum \limits_{\boldsymbol \delta} \frac{n_{\boldsymbol \delta}}{n}
\delta_j F'(i \boldsymbol \delta \cdot \boldsymbol \eta) = 2 \pi i
\frac{\alpha_N^{(j)}}{\sqrt{n}}, \quad j=1,\dots,\ell,
\end{equation}
where we wrote $F = - \log \hat{\rho}.$

\par Since $\max \limits_{\sigma,\sigma' \in X} |R(\sigma,\sigma')|$
is $\mathrm{P}$-almost surely of order $o(1)$ when $M=o(N)$ we can
apply Lemma 5.5 from \cite{BCMN2} and obtain that this system has a
unique solution
\begin{equation}
\boldsymbol \eta(\boldsymbol \alpha) = \frac{1}{2 \pi} B^{-1}
\boldsymbol \alpha \, \Bigl(1 + O\bigl(||\boldsymbol \alpha||_{2}^2
\bigr)\Bigr).
\end{equation}
Moreover, for sufficiently small~$\mu_1$,
\begin{align}\label{(a):ng:prob4}
& \iiint_{-\mu_1}^{\mu_1} e^{2 \pi n (i \mathbf{f} \cdot \boldsymbol
\alpha - \boldsymbol{\eta} \cdot \boldsymbol{\alpha}) } \prod
\limits_{\boldsymbol \delta} \hat{\rho} (\mathbf{f} \cdot
\boldsymbol \delta + i \boldsymbol \eta \cdot \boldsymbol \delta
)^{n_{\boldsymbol \delta}} \prod \limits_{j=1}^{\ell} df_j \nonumber \\
& \ \ = e^{-nG_{n,\ell}(\boldsymbol \alpha)} \biggl( \frac{1}{2 \pi
n}\biggr)^{\ell/2}    \bigl(1+O(n^{-1/2})+ O( a_N^2/n) +
O(R_{\max})\bigr),
\end{align}
where
\begin{equation}\label{(a):ng:G}
G_{n,\ell}(\boldsymbol \alpha) = \sum \limits_{\boldsymbol \delta}
\frac{n_{\boldsymbol \delta}}{n} F(i \boldsymbol \delta \cdot
\boldsymbol \eta(\boldsymbol \alpha)) + 2 \pi \boldsymbol
\eta(\boldsymbol \alpha) \cdot \boldsymbol \alpha.
\end{equation}

\noindent Therefore
\begin{align}\label{(ng):prob5}
\mathbb{P}&\bigl(H'_N(\sigma^j) \in (x_j, x_j + dx_j) \mbox{ for }
j=1,\dots,\ell \bigr) =  \biggl( \frac{b_N}{\sqrt{2
\pi}}\biggr)^{\ell} e^{-nG_{n,\ell}(\boldsymbol \alpha)} \prod
\limits_{j=1}^{\ell} dx_j \nonumber
\\ & \times  \bigl(1+O(n^{-1/2})+ O( a_N^2/n) +
O(R_{\max})\bigr) + O(b_N^\ell n^{\ell/2} e^{-\frac{1}{2} c_1
n_{\min}}).
\end{align}

\noindent Expanding $G_{n,\ell}$ we get the approximation
\begin{equation}\label{(ng):G_1}
n G_{n,\ell}(\boldsymbol \alpha) = \frac{n}{2} (\boldsymbol \alpha,
B^{-1}\boldsymbol \alpha) + O\Bigl(\frac{a_N^4}{n}\Bigr).
\end{equation}

\par By definition $a_N=O(\sqrt{M}).$ Under the assumptions of
Theorem~\ref{(ng):thm1} we get that $a_N =o(\sqrt[4]{N})$ for
$p=1$, that $a_N=o(\sqrt{N})$ for $p\geq 2$, and thus that
$a_N^4=o(n).$ It implies that asymptotically the joint density
\eqref{(ng):joint:den} is Gaussian. More precisely, it follows
from the equations \eqref{(ng):prob5} and \eqref{(ng):G_1} that
$\mathrm{P}$-a.s., for all collections
$(\sigma^1,\dots,\sigma^\ell) \in X^{\ell}$,
\begin{align}\label{(a):ng:prob6}
\mathbb{P}\bigl(H'_N(\sigma^j) & \in (x_j, x_j + dx_j) \mbox{ for }
j=1,\dots,\ell \bigr) =  \biggl( \frac{b_N}{\sqrt{2
\pi}}\biggr)^{\ell} e^{-n (\boldsymbol \alpha, B^{-1}\boldsymbol
\alpha)/2} \nonumber \\ & \times \prod \limits_{j=1}^{\ell} dx_j \,
\Bigl(1+O(R_{\max})+O\Bigl(\frac{a_N^4}{n}\Bigr)\Bigr) + O(b_N^\ell
n^{\ell/2} e^{-\frac{1}{2} c_1 n_{\min}}).
\end{align}
By Lemma~\ref{(ng):Cp:lem} the term $O(b_N^\ell n^{\ell/2}
e^{-\frac{1}{2} c_1 n_{\min}})$ is of order $o(e^{-cN^p})$ as $N \to
\infty.$ This finishes the proof of Theorem~\ref{(ng):thm1}.
\end{proof}

\medskip

\subsection{Breakdown of Universality}

In this subsection we follow the same strategy as we used in the
proof of Theorem \ref{main_thm_2} -- we fix a bounded set $A$ and
study the first three factorial moments of the random variable
$\mathcal{P}_N(A).$ In the case of the number partitioning problem
the following theorem implies that the Poisson convergence fails as
soon as $\limsup M/\sqrt{N} > 0.$
\begin{thm}[Number partitioning problem]
\label{(ng):thm3} Fix $p=1$ and let the Hamiltonian be given by
\eqref{(ng):Hamiltonian}. For every bounded Borel set $A$ we have
\begin{gather}
      \lim \mathbb{E} (\mathcal{P}_N(A))_1 =
      \mu(A)e^{-4 c_4 \varepsilon^2 \log^2 2}, \quad \mathrm{P}\mbox{-a.s.}
\end{gather}
Moreover, if $\ \limsup \frac{M(N)}{N} = \varepsilon < \infty$
then $\mathrm{P}$-a.s.

$(i)$ $\limsup \mathbb{E} (\mathcal{P}_N(A))_2 = e^{2 \varepsilon^2
\log 2 - 32 c_4 \varepsilon^2 \log^2 2} (\mu(A))^2;$

$(ii)$ $\limsup \mathbb{E} (\mathcal{P}_N(A))_3 < \infty.$
\end{thm}
Therefore the limit of the ratio of the second factorial moment to
the square of the first is
\begin{equation}
\frac{\mathbb{E}(\mathcal{P}_N(A))_2}{\mathbb{E}(\mathcal{P}_N(A))^2_1}
\to e^{2 \varepsilon^2 \log^2 2 - 24 c_4 \varepsilon^2 \log^2 2 } =
e^{2 \varepsilon^2 \log^2 2(1 - 12 c_4)}.
\end{equation}
Taking into account that $c_4 < \frac{1}{12}$ we conclude that the
ratio is strictly larger than one and thus there is no Poisson
convergence for $\varepsilon >0.$

\medskip

And in the case of the Sherrington-Kirkpatrick model the failure of
Poisson convergence follows from
\begin{thm}[Sherrington-Kirkpatrick
model]\label{(ng):thm2} Fix $p=2$ and let the Hamiltonian be given
by \eqref{(ng):Hamiltonian}. For every bounded Borel set $A$
\begin{gather}
\label{1st'}
      \lim \mathbb{E} (\mathcal{P}_N(A))_1 =
      \mu(A)e^{-4 c_4 \varepsilon^2 \log^2 2}, \quad \mathrm{P}\mbox{-a.s.}
\end{gather}
Moreover, if $\ \limsup \frac{M(N)}{N} = \varepsilon < \frac{1}{8
\log 2}$ then $\mathrm{P}$-a.s.

(i) $\limsup \mathbb{E} (\mathcal{P}_N(A))_2 = \frac{e^{- 32 c_4
\varepsilon^2 \log^2 2}}{\sqrt{1 - 4 \varepsilon \log
2}}(\mu(A))^2;$

(ii) $\limsup \mathbb{E} (\mathcal{P}_N(A))_3 < \infty.$
\end{thm}
The ratio of the second factorial moment to the square of the first
moment is
\begin{equation}
\frac{\mathbb{E}(\mathcal{P}_N(A))_2}{\mathbb{E}(\mathcal{P}_N(A))^2_1}
\to  \frac{e^{- 24 c_4 \varepsilon^2 \log^2 2}}{\sqrt{1 - 4
\varepsilon \log 2}}.
\end{equation}
For $\varepsilon> 0$ the above ratio is strictly larger than one and
thus convergence to a Poisson point process fails.

We will give the proof of Theorem \ref{(ng):thm2} only since the
case $p=1$ is based on essentially the same computations.

\begin{proof}[Proof of Theorem \ref{(ng):thm2}]
As in the proof of Theorem~\ref{(c):thm1} we successively prove the
statement on the first, second, and third moment. To simplify our
computations we will assume that $M = \varepsilon N$ (Using the
monotonicity argument, the case of general sequences $M(N)$ can be
analyzed just as in Theorem \ref{(c):thm1} of Section \ref{(c)}.).


\bigskip
\noindent 1.  {\bf  First moment estimate.}

Following the same steps as in Subsection \ref{(ng.a)} we
approximate the density of $H'_N(\sigma)$ by
\begin{equation}
e^{-nG_{n,1}(\alpha_N)} \frac{1}{\sqrt{2 \pi}} \Bigl(1 +
O\Bigl(\frac{1}{\sqrt{N}}\Bigr) + O\Bigl(\frac{a_N^2}{n}\Bigr)
\Bigr) + O(e^{-c_1 n}),
\end{equation}
where according to the notation introduced above $n=N^2.$

\par In the case $M = \varepsilon N$ we need  a more
precise approximation of the function $G_{n,\ell}$ than given by
formula \eqref{(ng):G_1}. Expanding the solution of the system
\eqref{(a):ng:system:G} as
\begin{equation}\label{(b):ng:eta}
\boldsymbol \eta(\boldsymbol \alpha) = \frac{1}{2 \pi} B^{-1}
\boldsymbol \alpha + \frac{4 c_4}{2 \pi} \sum \limits_{\boldsymbol
\delta} \frac{n_{\boldsymbol \delta}}{n}(\boldsymbol \delta, B^{-1}
\boldsymbol \alpha)^3 B^{-1} \boldsymbol \delta + O(||\alpha||^5)
\end{equation}
and applying \eqref{(b):ng:eta} we obtain from \eqref{(a):ng:G} and
\eqref{(ng):F} that
\begin{align}\label{(b):ng:G}
G_{n,\ell}(\boldsymbol \alpha) & = -\frac{(2 \pi)^2}{2} \sum
\limits_{\boldsymbol \delta} \frac{n_{\boldsymbol
\delta}}{n}(\boldsymbol \delta \cdot \boldsymbol \eta)^2 \nonumber
\\ & + c_4 (2\pi)^4 \sum \limits_{\boldsymbol \delta}
\frac{n_{\boldsymbol \delta}}{n} (\boldsymbol \delta \cdot
\boldsymbol \eta)^4 + 2 \pi (\boldsymbol \eta \cdot \boldsymbol
\alpha) + O(||\boldsymbol \alpha||^6)
\nonumber \\
& = \frac{1}{2}(\boldsymbol \alpha, B^{-1} \boldsymbol \alpha) + c_4
\sum \limits_{\boldsymbol \delta} \frac{n_{\boldsymbol \delta}}{n}
(\boldsymbol \delta, B^{-1} \boldsymbol \alpha)^4 + O(||\boldsymbol
\alpha||^6).
\end{align}

\noindent Using the approximation for $G_{n,1}$ given by formula
\eqref{(b):ng:G}, we obtain
\begin{equation}
nG_{n,1} = \frac{\alpha_N^2}{2}  + c_4 \frac{\alpha_N^4}{n} +
O\Bigl(\frac{a_N^6}{n^2}\Bigr),
\end{equation}
where $\alpha_N = a_N + b_N \, x.$ Since $\alpha_N^4/n = 4 c_4
\varepsilon^2 \log^2 2 (1+o(1)),$ we see that up to an
$\varepsilon$-dependent multiplier the density of $H'_N(\sigma)$ is
given by the normal density, more precisely, the density of
$H'_N(\sigma)$ is
\begin{align}
\frac{1}{\sqrt{2 \pi}} e^{-\alpha_N^2/2} \, e^{-4 c_4 \varepsilon^2
\log^2 2} \Bigl(1 + O\Bigl(\frac{1}{\sqrt{N}}\Bigr) +
O\Bigl(\frac{\alpha_N^2}{n}\Bigr) \Bigr) + O(e^{-c_1 n}).
\end{align}
Therefore the first factorial moment of $\mathcal{P}_N(A)$ is
\begin{equation}\label{(c):ng:first_mom}
\mu(A)e^{-4 c_4 \varepsilon^2 \log^2 2}(1+o(1)) + o(1)
\end{equation}
and \eqref{1st'} is proven.

\bigskip
\noindent 2.  {\bf  Second moment estimate.} When analyzing the
second moment
\begin{equation}\label{(ng):c:moment}
\mathbb{E}(\mathcal{P}_N(A))_2 = \sum \limits_{
(\sigma^{1},\sigma^{2}) \in X^{2}} \mathbb{P}\left(
H'_{N}(\sigma^{1}) \in A,H'_{N}(\sigma^{2}) \in A\right ),
\end{equation}
it is useful to distinguish  between ``typical'' and ``atypical''
sets of configurations $(\sigma^1,\dots,\sigma^\ell) \in X^{\ell},$
a notion introduced in \cite{BCMN1}.

\par Take some sequence $\theta_N \to 0$ such that $N \theta_N^2
\to \infty$ and consider the $\ell \times n$ matrix
$C_p(\sigma^1,\dots,\sigma^\ell),$ introduced in Subsection
\ref{(ng.a)}. Then all but a vanishing fraction of the
configurations $(\sigma^1,\dots,\sigma^\ell) \in S_N^{\ell}$ obey
the condition
\begin{equation}\label{(ng):typical}
\max \limits_{\boldsymbol \delta \in \{-1,1\}^{\ell} } \Bigl| \,
n_{\boldsymbol{\delta}} - \frac{n}{2^{\ell}} \Bigr| \leq n \theta_N.
\end{equation}

\par When $M=o(N)$, Lemma \ref{(ng):Cp:lem} guarantees that for a
properly chosen sequences $\theta_N$, $\mathrm{P}$-almost surely,
all the sampled sets $(\sigma^1,\dots,\sigma^\ell) \in X^{\ell}$
obey condition \eqref{(ng):typical}. It is no longer the case when
$M = \varepsilon N$ and thus we have to consider the contribution
from the sets violating \eqref{(ng):typical}.

\par Fix $p=2$ and $\ell = 2.$ Let a sequence $\theta_N \to 0$ be given, and define
\begin{align}\label{(c):ng:I}
I = &\sum \limits_{\sigma^1,\sigma^2} \mathbb{P}\bigl(
H'_N(\sigma^1) \in A, H'_N(\sigma^2) \in A \bigr),
\end{align}
where the sum runs over all pairs of distinct configurations
$(\sigma^1,\sigma^2) \in X^{2}$ satisfying condition
\eqref{(ng):typical} (the so called ``typical" configurations). Also
define
\begin{align}\label{(c):ng:II}
II = &\sum \limits_{\sigma^1,\sigma^2} \mathbb{P}\bigl(
H'_N(\sigma^1) \in A, H'_N(\sigma^2) \in A \bigr),
\end{align}
where the sum is over all pairs of distinct configurations
$(\sigma^1,\sigma^2) \in X^{2}$ violating \eqref{(ng):typical} (the
``atypical" configurations). For later use we introduce the
quantities $I_g$ and $II_g$ -- the analogs of the variables $I$ and
$II$ in the case where the random variables $(g_{i_1,i_2})_{1 \leq
i_1,i_2 \leq N}$ are i.i.d. standard normals.
\medskip

\begin{lem}\label{(c):ng:I&II}
Let $\theta_N = \frac{1}{\sqrt{N} \log^2 N}.$ If $\varepsilon \in
(0,\frac{1}{4 \log 2})$ then $\mathrm{P}$-almost surely $II = o(1)$
and
\begin{equation}
I = \frac{e^{ - 32 c_4 \varepsilon^2 \log^2 2}
(\mu(A))^2}{\sqrt{1 - 4 \varepsilon \log 2}}(1+o(1)).
\end{equation}
\end{lem}

\begin{proof}

\par To prove that $II$ is $o(1)$ let us
bound the quantity $II$ by $II_g$ and show that $II_g$ is almost
surely negligible for $\theta_N=\frac{1}{\sqrt{N} \log^2 N}.$ We
start with the proof of the second statement, for which we will need
the following simple observation.
\par Consider a sequence $\lambda_N \to
0$ such that $N \lambda_N \to \infty$ and define the set of
configurations $\sigma \in \{-1,1\}^N$ with almost equal number of
spins equal to $1$ and to~$-1:$
\begin{equation}
\mathcal{T}_N = \Bigl\{ \sigma \in S_N: \bigl| \# \{\sigma_i=1\} -
\# \{\sigma_i=-1\} \bigr| \leq N \lambda_N \Bigr\}.
\end{equation}
It is not hard to prove that configurations $\sigma^1, \sigma^2 \in
\mathcal{T}_N$ with overlap $|R(\sigma^1,\sigma^2)| \leq \lambda_N$
must satisfy \eqref{(ng):typical} with $\theta_N = \lambda_N^2.$
Therefore the set of pairs $(\sigma^1, \sigma^2) \in X^{2}$
violating condition \eqref{(ng):typical} with $\theta_N =
\lambda_N^2$ is contained in the set
\begin{equation}\label{(ng):atyp}
\Bigl\{ (\sigma^1, \sigma^2) \in X^{2}: |R(\sigma^1,\sigma^2)|
> \lambda_N\mbox{ or }\sigma^1 \in \mathcal{T}_N^c \mbox{ or }\sigma^2
> \in \mathcal{T}_N^c \Bigr\}.
\end{equation}

\noindent Thus to prove that $II_{g}$ is $o(1)$ for $\theta_N =
\frac{1}{\sqrt{N} \log^2 N}$ it suffices to prove that
\begin{equation}\label{(ng):sum_atyp}
\sum \limits_{\sigma^1,\sigma^2} \mathbb{P}(H'_N(\sigma^1) \in A,
H'_N(\sigma^2) \in A) = o(1),
\end{equation}
where the summation is over all pairs of distinct configurations
contained in the set \eqref{(ng):atyp} with $\lambda_N =
\sqrt{\theta_N} = \frac{1}{N^{1/4} \log N}.$

\par Let us prove \eqref{(ng):sum_atyp}. Since we already proved in
Section \ref{(c)} that the contribution from the set
\begin{equation}
\bigl\{(\sigma^1,\sigma^2) \in X^{2}: |R(\sigma^1,\sigma^2)|
> \lambda_N \bigr\}
\end{equation}
to the sum \eqref{(ng):c:moment} is negligible, it is enough to
consider the sum \eqref{(ng):c:moment} restricted to the set

\noindent
\begin{equation}\label{(ng):sum_atyp1}
\Bigl\{(\sigma^1,\sigma^2) \in X^{2}: \sigma^1 \mbox{ or } \sigma^2
\in (\mathcal{T}_N^c)^X \mbox{ and } |R(\sigma^1, \sigma^2)| <
\lambda_N \Bigr\}.
\end{equation}

\noindent By Stirling's formula we obtain that $|\mathcal{T}_N^c| =
\sqrt{\frac{2}{\pi N}}2^{N}
e^{-N\mathcal{J}(\lambda_N)}(1+O(\lambda_N^2))$ and using this fact
one can prove, proceeding as in part $(i)$ of Theorem
\ref{(ce):thm3}, that
\begin{equation}
|(\mathcal{T}_N^c)^X| = \mathrm{E} |(\mathcal{T}_N^c)^X| (1+o(1)) =
\sqrt{\frac{2}{\pi N}} \, 2^{M}
e^{-N\mathcal{J}(\lambda_N)}(1+o(1)).
\end{equation}
Thus for large enough $N$ the contribution from the set
\eqref{(ng):sum_atyp1} to the sum \eqref{(ng):c:moment} is bounded
by
\begin{equation}\label{(ng):sum_atyp2}
2 \sqrt{\frac{2}{\pi N}}2^{2 M} e^{-N\mathcal{J}(\lambda_N)}
\frac{C}{2^{2M}} e^{\frac{2b{12}}{1+b_{12}}(M \log 2 - 1/2 \log M)},
\end{equation}
where the constant $C$ is from \eqref{(c):g:int2}. Using that $N
\mathcal{J}( \lambda_N )=\frac{1}{2}N\lambda_N^2 +
O\Bigl(\frac{1}{\log^4 N}\Bigr)$ we can further bound
\eqref{(ng):sum_atyp2} by
\begin{equation}
2 C \sqrt{\frac{2}{\pi N}} e^{-\frac{1}{2}N \lambda_N^2(1 - 4
\varepsilon \log 2) },
\end{equation}
which is $o(1)$ by the choice of $\lambda_N$ and $\varepsilon.$

\par Our next step is to bound the sum $II$ by $II_g.$
For this purpose we need to give an estimate of the joint density of
$H'_N(\sigma^1),H'_N(\sigma^2)$ that would be valid also for pairs
$(\sigma^1,\sigma^2)$ violating condition \eqref{(ng):typical}. As
we already noted for such $(\sigma^1,\sigma^2)$ the results of
Subsection \ref{(ng.a)} cannot be applied directly since it is no
longer true that $\max \limits_{\sigma,\sigma' \in X}
|R(\sigma,\sigma')|$ is $o(1).$ Fortunately, we have only $2\times
2$ covariance matrix $B(\sigma^1,\sigma^2)$ and using this fact
we can easily adapt the results of Subsection \ref{(ng.a)} to the
case where $\max \limits_{\sigma,\sigma' \in X} |R(\sigma,\sigma')|$
is not $o(1).$ We start with formula \eqref{(ng):prob} which, in the
case $\ell = 2$, can be rewritten as
\begin{align}
\mathbb{P}\bigl( H'_N(&\sigma^j) \in (x_j, x_j + dx_j) \mbox{ for }
j=1,2 \bigr) \nonumber \\ & = b_N^2 n \, dx_1 dx_2
\iint_{-\infty}^{\infty} \prod \limits_{1 \le i_1, i_2 \le N}
\hat{\rho}(v_{i_1,i_2}) e^{2 \pi i \sqrt{n} \bigl( f_1
\alpha_N^{(1)} + f_2 \alpha_N^{(2)}\bigr)} df_1 df_2.
\end{align}

\noindent We can rewrite the integral in the above expression as
\begin{align}\label{(ng):int}
\iint_{-\infty}^{\infty} \prod \limits_{ \boldsymbol{\delta}}
\hat{\rho}({\mathbf f}\cdot \boldsymbol{\delta}
)^{n_{\boldsymbol{\delta}}} \, e^{2 \pi i n \mathbf{f} \cdot
\boldsymbol{\alpha}} df_1 df_2,
\end{align}
where $\boldsymbol \delta \in \{-1,1\}^2.$ Since the function
$\hat{\rho}$ is even we obtain
\begin{equation}
\prod \limits_{ \boldsymbol{\delta}} \hat{\rho}({\mathbf f}\cdot
\boldsymbol{\delta} )^{n_{\boldsymbol{\delta}}} =
\hat{\rho}(f_1+f_2)^{n_{(1,1)} +
n_{(-1,-1)}}\hat{\rho}(f_1-f_2)^{n_{(1,-1)} + n_{(-1,1)}}.
\end{equation}
One obvious relation between $n_{(1,1)},n_{(-1,-1)},n_{(1,-1)}$ and
$n_{(-1,1)}$ is
\begin{equation}
n_{(1,1)} + n_{(-1,-1)}+ n_{(1,-1)}+n_{(-1,1)} = n.
\end{equation}
The other one we obtain by noting that
\begin{align}
n_{(1,1)} + & n_{(-1,-1)} - n_{(-1,1)} - n_{(1,-1)} = \sum
\limits_{i,j} \sigma_i^1 \sigma_j^1 \sigma_i^2 \sigma_j^2 =
nR_{12}^2.
\end{align}

\noindent Therefore
\begin{equation}
\left\{
\begin{array}{lr}
n_{(1,1)}+n_{(-1,-1)} = \frac{1}{2}n(1+R_{12}^2),\\
n_{(1,-1)}+n_{(-1,1)} = \frac{1}{2}n(1-R_{12}^2 ).
\end{array}
\right.
\end{equation}
By Theorem \ref{(ce):thm3} we conclude that $\mathrm{P}$-a.s. $\max
\limits_{\sigma,\sigma' \in X} |R(\sigma,\sigma')| < 1$ and
therefore, for some positive constant $c$,
\begin{equation}\label{(c):ng:n_min}
n_{(1,1)}+n_{(-1,-1)} \geq n_{(1,-1)}+n_{(-1,1)} \geq cn.
\end{equation}
The above inequality allows us to approximate \eqref{(ng):int} by
\begin{align}
\iint_{-\mu_1}^{\mu_1} e^{2 \pi n (i \mathbf{f} \cdot \boldsymbol
\alpha - \boldsymbol{\eta} \cdot \boldsymbol{\alpha}) } \prod
\limits_{\boldsymbol \delta} \hat{\rho} (\mathbf{f} \cdot
\boldsymbol \delta + i \boldsymbol \eta \cdot \boldsymbol \delta
)^{n_{\boldsymbol \delta}} \, df_1 df_2 + O\bigl(e^{-c_1 n}\bigr).
\end{align}

\noindent Adapting the proof of Lemma 5.5 from \cite{BCMN2} we get
\begin{align}
\iint_{-\mu_1}^{\mu_1} & e^{2 \pi n (i \mathbf{f} \cdot \boldsymbol
\alpha - \boldsymbol{\eta} \cdot \boldsymbol{\alpha}) } \prod
\limits_{\boldsymbol \delta} \hat{\rho} (\mathbf{f} \cdot
\boldsymbol \delta + i \boldsymbol \eta \cdot \boldsymbol \delta
)^{n_{\boldsymbol \delta}} \, df_1 df_2 \nonumber \\
& = e^{-nG_{n,2}(\boldsymbol
\alpha)}\frac{\sqrt{\det{B}(\sigma^1,\sigma^2)}}{2 \pi n} \Bigl(1 +
O\Bigl(\frac{1}{\sqrt{n}}\Bigr) +
O\Bigl(\frac{a_N^2}{n}\Bigr)\Bigr).
\end{align}
Finally, we obtain that
\begin{align}
\mathbb{P}\bigl(H'_N(&\sigma^j)\in (x_j,x_j+dx_j) \mbox{ for }
j=1,2 \bigr) = dx_1 dx_2 \, b_n e^{-n G_{n,2}(\boldsymbol \alpha)}
\nonumber
\\ & \times \frac{\sqrt{\det{B}}}{2 \pi }
\Bigl(1 + O\Bigl(\frac{1}{\sqrt{n}}\Bigr) +
O\Bigl(\frac{a_N^2}{n}\Bigr)\Bigr) + O\bigl(b_N^2 n e^{-c_1n}\bigr).
\end{align}

\par From \eqref{(b):ng:G} we see that for some constant $C$
\begin{align}
n G_{n,2} & = \frac{n}{2}(\boldsymbol \alpha, B^{-1} \boldsymbol
\alpha) + \frac{16 c_4}{(1+R_{12}^2)^3} \frac{\alpha_N^4}{n} +
O\Bigl(\frac{a_N^6}{n^2}\Bigr) \ge \frac{n}{2}(\boldsymbol \alpha,
B^{-1} \boldsymbol \alpha) + C,
\end{align}
and thus the joint density of $H'_N(\sigma^1),H'_N(\sigma^2)$ is
bounded by
\begin{equation}\label{(c):ng:fin_dens}
\frac{b_N^2 \sqrt{\det{B}}}{2 \pi} \, e^{-n(\boldsymbol\alpha,
B^{-1} \boldsymbol \alpha)/2 - C}\Bigl(1 +
O\Bigl(\frac{1}{\sqrt{n}}\Bigr) + O\Bigl(\frac{a_N^2}{n}\Bigr)\Bigr)
+ O\bigl(b_N^2 n e^{-c_1 n}\bigr).
\end{equation}

\noindent This last bound for the joint density clearly implies that
the sum $II$ could be bounded by $II_{g}$ plus an error resulting
from the second term. But the cumulative error coming from the
second term is of order $O\bigl(2^{2M}b_N^2 n e^{-c_1 n}\bigr)$
which is negligible even in the case $\limsup M/N  >0.$
\smallskip

\par To prove the second statement of the lemma we will approximate
the sum $I$ by $I_g,$ which was already calculated in Section
\ref{(c)}. We first notice that for sequences
$(\sigma^1,\dots,\sigma^\ell)$ satisfying condition
\eqref{(ng):typical} $R_{\max}(\sigma^1,\dots,\sigma^\ell ) =
O(\theta_N).$ Furthermore, for configurations
$(\sigma^1,\dots,\sigma^\ell)$ obeying condition
\eqref{(ng):typical} it is possible to derive from \eqref{(b):ng:G}
that
\begin{equation}\label{(ng):c:G}
G_{n,\ell}(\boldsymbol \alpha) = \frac{1}{2}(\boldsymbol \alpha,
B^{-1} \boldsymbol \alpha) + c_4 \ell (1+3(\ell -1))\frac{a_N^4}{n}
+ O\Bigl(\frac{a_N^6}{n^3}\Bigr) + O\Bigl(\frac{a_N^2}{n^2}
\theta_N\Bigr).
\end{equation}
For the details of the derivation we refer to Subsection 5.4 of
\cite{BCMN2} and in particular to formula (5.57) in there. Using
formula \eqref{(ng):c:G} with $\ell = 2$ and substituting it into
\eqref{(ng):prob5} we obtain that
\begin{equation}
I = I_g \, e^{- 32 c_4 \varepsilon^2 \log^2 2} \Bigl(1 +
O\Bigl(\frac{1}{\sqrt{n}}\Bigr) +
O\Bigl(\frac{\alpha_N^2}{n}\Bigr)\Bigr).
\end{equation}
This finishes the proof of Lemma \ref{(c):ng:I&II}.
\end{proof}

To conclude the calculation of the second moment we notice that
summing $I$ and $II$ we get the second factorial moment
\begin{equation}\label{(c):ng:second_mom}
\mathbb{E}(\mathcal{P}_N(A))_2 = \frac{e^{- 32 c_4 \varepsilon^2
\log^2 2} (\mu(A))^2}{\sqrt{1 - 4 \varepsilon \log 2}}(1+o(1)).
\end{equation}
Assertion (i) of Theorem \ref{(ng):thm2} is thus proven.


\bigskip
\noindent 3.  {\bf  Third moment estimate.} To deal with the third
moment
\begin{equation}\label{(ng):c:moment:3}
\mathbb{E}(\mathcal{P}_N(A))_3 = \sum \limits_{
(\sigma^{1},\sigma^{2},\sigma^{3}) \in X^{3}} \mathbb{P}\left(
H'_{N}(\sigma^{1}) \in A,H'_{N}(\sigma^{2}) \in A,
H'_{N}(\sigma^{3}) \in A\right ),
\end{equation}
we use the same strategy as we used to calculate the second moment.
In particular, we fix $\ell = 3$ and split the sum
\eqref{(ng):c:moment:3} in two parts:
\begin{align}
I = &\sum \limits_{\sigma^1,\sigma^2,\sigma^3} \mathbb{P}\bigl(
H'_N(\sigma^1) \in A, H'_N(\sigma^2) \in A, H'_N(\sigma^3) \in A
\bigr),
\end{align}
where the sum runs over all sequences of distinct configurations
$(\sigma^1,\sigma^2,\sigma^3) \in X^{3}$ satisfying condition
\eqref{(ng):typical}, and
\begin{align}
II = &\sum \limits_{\sigma^1,\sigma^2,\sigma^3} \mathbb{P}\bigl(
H'_N(\sigma^1) \in A, H'_N(\sigma^2) \in A, H'_N(\sigma^3) \in A
\bigr),
\end{align}
where the sum is over all sequences of distinct configurations
$(\sigma^1,\sigma^2,\sigma^3) \in X^{3}$ violating
\eqref{(ng):typical}.

By exactly the same argument as in the calculation of the second
moment the contribution from the ``typical'' collections, $I,$ is
bounded. We therefore concentrate on the analysis of the
contribution from the ``atypical'' collections, $II.$ Since we are
interested only in the estimate of the third moment from above it
suffices to bound the joint density $\mathbb{P}\bigl(H'_N(\sigma^j)
\in (x_j, x_j + dx_j) \mbox{ for } j=1,2,3 \bigr).$ We start with
formula \eqref{(ng):prob} which, in the case $\ell = 3$, can be
rewritten as
\begin{align}\label{(ng):int:3}
\mathbb{P}\bigl( H'_N(\sigma^j) & \in (x_j, x_j + dx_j) \mbox{ for }
j=1,2,3 \bigr) = b_N^3 n \, dx_1 dx_2 dx_3 \nonumber \\ & \times
\iiint_{-\infty}^{\infty} \prod \limits_{ i_1, i_2, i_3}
\hat{\rho}(v_{i_1,i_2,i_3}) e^{2 \pi i \sqrt{n}\bigl( f_1
\alpha_N^{(1)} + f_2 \alpha_N^{(2)} + f_3 \alpha_N^{(3)} \bigr)}
df_1 df_2 df_3.
\end{align}

\noindent We can rewrite the integral in the above expression as
\begin{align}
\iiint_{-\infty}^{\infty} \prod \limits_{ \boldsymbol{\delta}}
\hat{\rho}({\mathbf f}\cdot \boldsymbol{\delta}
)^{n_{\boldsymbol{\delta}}} \, e^{2 \pi i n \mathbf{f} \cdot
\boldsymbol{\alpha}} df_1 df_2 df_3,
\end{align}
where $\boldsymbol \delta \in \{-1,1\}^3.$ Since the function
$\hat{\rho}$ is even we obtain
\begin{align}
\prod \limits_{ \boldsymbol{\delta}} \hat{\rho}({\mathbf f}\cdot
\boldsymbol{\delta} )^{n_{\boldsymbol{\delta}}} & =
\hat{\rho}(f_1+f_2+f_3)^{n_1}\hat{\rho}(f_1+f_2-f_3)^{n_2} \nonumber \\
& \times \hat{\rho}(f_1-f_2+f_3)^{n_3}
\hat{\rho}(f_1-f_2-f_3)^{n_4},
\end{align}
where
\begin{equation}
\left\{
\begin{array}{lr}
n_1 = n_{(1,1,1)} + n_{(-1,-1,-1)},\\
n_2 = n_{(1,1,-1)} + n_{(-1,-1,1)},\\
n_3 = n_{(1,-1,1)} + n_{(-1,1,-1)},\\
n_4 = n_{(-1,1,1)} + n_{(1,-1,-1)}.
\end{array}
\right.
\end{equation}
By definition of the matrix $C_2(\sigma^1,\sigma^2,\sigma^3)$ we
have
\begin{equation}
\left\{
\begin{array}{lr}
n_1 + n_2 - n_3 - n_4 = nR_{12}^2,\\
n_1 - n_2 - n_3 + n_4 = nR_{23}^2,\\
n_1 - n_2 + n_3 - n_4 = nR_{31}^2,\\
n_1 + n_2 + n_3 + n_4 = n.\\
\end{array}
\right.
\end{equation}
Solving the system
\begin{equation}
\left\{
\begin{array}{lr}
n_1 = \frac{1}{4} n ( 1 + R_{12}^2 + R_{23}^2 + R_{31}^2),\\
n_2 = \frac{1}{4} n ( 1 + R_{12}^2 - R_{23}^2 - R_{31}^2),\\
n_3 = \frac{1}{4} n ( 1 - R_{12}^2 - R_{23}^2 + R_{31}^2),\\
n_4 = \frac{1}{4} n ( 1 - R_{12}^2 + R_{23}^2 - R_{31}^2).\\
\end{array}
\right.
\end{equation}
From Theorem \ref{(ce):thm3} we obtain that $\mathrm{P}$-almost
surely
\begin{equation}\label{(ng):3:moment:bound:R}
\limsup \limits_{N \to \infty} \max \limits_{\sigma^1,\sigma^2 \in
X} \mathcal{J}\bigl(R(\sigma^1,\sigma^2)\bigr) \le \varepsilon \log
2.
\end{equation}
Since the function $\mathcal{J}$ is monotone we obtain from
\eqref{(ng):3:moment:bound:R} and from assumption $\varepsilon <
\frac{1}{8 \log 2}$ that $\mathrm{P}$-a.s.
\begin{equation}
\limsup \limits_{N \to \infty }\max \limits_{\sigma^1,\sigma^2 \in
X} |R(\sigma^1,\sigma^2)| < \frac{1}{2}.
\end{equation}
It implies that $\min \{ n_1,n_2,n_3,n_4 \} \ge cn$ for some
positive constant $c.$ It allows us to approximate the integral in
\eqref{(ng):int:3} by
\begin{align}
\iiint_{-\mu_1}^{\mu_1} e^{2 \pi n (i \mathbf{f} \cdot \boldsymbol
\alpha - \boldsymbol{\eta} \cdot \boldsymbol{\alpha}) } \prod
\limits_{\boldsymbol \delta} \hat{\rho} (\mathbf{f} \cdot
\boldsymbol \delta + i \boldsymbol \eta \cdot \boldsymbol \delta
)^{n_{\boldsymbol \delta}} \, df_1 df_2 df_3 + O\bigl(e^{-c
n}\bigr).
\end{align}
Adapting the proof of Lemma 5.5 from \cite{BCMN2} we get
\begin{align}
\iiint_{-\mu_1}^{\mu_1} & e^{2 \pi n (i \mathbf{f} \cdot \boldsymbol
\alpha - \boldsymbol{\eta} \cdot \boldsymbol{\alpha}) } \prod
\limits_{\boldsymbol \delta} \hat{\rho} (\mathbf{f} \cdot
\boldsymbol \delta + i \boldsymbol \eta \cdot \boldsymbol \delta
)^{n_{\boldsymbol \delta}} \, df_1 df_2 df_3 \nonumber \\
& = e^{-nG_{n,3}(\boldsymbol
\alpha)}\frac{\sqrt{\det{B}(\sigma^1,\sigma^2)}}{(2 \pi n)^{3/2}}
\Bigl(1 + O\Bigl(\frac{1}{\sqrt{n}}\Bigr) +
O\Bigl(\frac{a_N^2}{n}\Bigr)\Bigr).
\end{align}
Next, after a little algebra, one can derive from \eqref{(b):ng:G}
that for some constant $C$
\begin{equation}
n G_{n,3} \ge \frac{n}{2}(\boldsymbol \alpha, B^{-1} \boldsymbol
\alpha) + C
\end{equation}
and this estimate is enough to bound the joint density of
$H'_N(\sigma^1),H'_N(\sigma^2),H'_N(\sigma^3)$ by the joint density
of $Z_N(\sigma^1), Z_N(\sigma^2), Z_N(\sigma^3).$ Thus Theorem
\ref{(ng):thm2} is proved.
\end{proof}

\bigskip

\subsection*{Acknowledgements}  We thank J.Cerny for a helpful reading of the first
versions of the manuscript. V. Gayrard thanks the Chair of
Stochastic Modeling of the \'Ecole Polytechnique F\'ed\'erale of
Lausanne for financial support.

\end{document}